%% file: root_v3.tex
\documentclass[11pt]{article}
\usepackage[margin=1in]{geometry}
\usepackage{bbm}
\usepackage[usenames]{color}
\usepackage{enumerate}
\usepackage{url}
\usepackage{subfigure}
\usepackage{amsfonts,mathrsfs}
\usepackage{amssymb,amsmath}
\usepackage{verbatim}
\usepackage{acronym}
\usepackage{mathtools}
\DeclarePairedDelimiter{\ceil}{\lceil}{\rceil}
\usepackage{graphicx}
\usepackage{todonotes}

\usepackage{enumerate}
\usepackage{amsthm}

\newcommand{\remove}[1]{}

\def\beh#1{{\color{black}#1}}
\def\roh#1{{{\color{blue}#1}}}
\def\fskip#1{}
\def\R{\mathbb{R}}

\def\allzero{\mathbf{0}}
\def\allone{\mathbf{1}}

\usepackage{amsmath}
\DeclareMathOperator*{\argmax}{arg\,max}
\DeclareMathOperator*{\argmin}{arg\,min}

\newtheorem{theorem}{Theorem}

\newtheorem{corollary}{Corollary}

\newtheorem{definition}{Definition}

\newtheorem{lemma}{Lemma}

\newtheorem{proposition}[theorem]{Proposition}
\newtheorem{remark}{Remark}

\renewcommand{\theenumi}{\arabic{enumi}}
\newcommand{\norm}[1]{\left\lVert#1\right\rVert}

\usepackage[normalem]{ulem}
\usepackage{caption}
\usepackage{ragged2e}
\usepackage{threeparttable}

\def\Gt{G}
\def\Cgph{\mathcal C(G_{ph})}
\def\At{A}
\def\Dt{D}

\def\Nt{N}
\def\Gph{G_{ph}}
\def\Eph{E_{ph}}
\def\Gn{\mathcal{G}_n}
\def\Pr{\text{Pr}}

\def\spann{\text{span}}

\begin{document}

\title{On the Convergence Properties of Social Hegselmann-Krause Dynamics
\author{Rohit Parasnis, Massimo Franceschetti, and Behrouz Touri\\
  Department of Electrical and Computer Engineering\\
     University of California San Diego, CA 92093\\
     Email: \{rparasni, mfranceschetti, btouri\}@eng.ucsd.edu}}
     \date{}
\maketitle

\begin{abstract}
We \beh{study} the convergence properties of \textit{Social Hegselmann-Krause dynamics}, a {variant} of the Hegselmann-Krause (HK) model of opinion dynamics where a physical connectivity graph that accounts for the extrinsic factors that could prevent interaction between certain pairs of agents is incorporated. \beh{ As opposed to the original HK dynamics (which terminate in finite time), we show that for any underlying connected and incomplete graph}, under a certain mild assumption, the expected termination time of social HK dynamics is infinity. We then investigate the rate of convergence to \beh{the} steady state, and provide bounds on the maximum $\epsilon$-convergence time in terms of the properties of the physical connectivity graph. We extend this discussion and observe that for almost all $n$, there exists an $n$-vertex physical connectivity graph on which social HK dynamics may not even $\epsilon$-converge to the steady state within a bounded time frame. We then provide nearly tight necessary and sufficient conditions for arbitrarily slow merging (a phenomenon that is essential for arbitrarily slow $\epsilon$-convergence to the steady state). 
Using the necessary conditions, we show that complete $r$-partite graphs have bounded $\epsilon$-convergence times. 

\end{abstract}


\input{introduction}
\input{formulation}
\input{terminationtime}
\input{upperbound}
\input{arbit_slow_conv_2}
\input{cdc_part}
\input{conclusion}
\bibliographystyle{abbrv}
\bibliography{bib}
\end{document}

%% file: introduction.tex
\section{INTRODUCTION}

With social networks gaining omnipresence and their associated datasets becoming accessible to the public, opinion dynamics has attracted researchers from a range of disciplines in recent times~\cite{etesami2013termination}. Besides having social scientific applications such as forecasting election results~\cite{sobkowicz2016quantitative}, opinion dynamics models are also used in engineering problems such as \textit{distributed rendezvous in a robotic network}~\cite{martinez2007synchronous}.

Among the existing models, confidence-based models form a noteworthy class. In particular, a well-known bounded-confidence model proposed in~\cite{hegselmann2002opinion}, also known as the Hegselmann-Krause model (referred as the \textit{HK model} from here on), has garnered a lot of interest in the last two decades. Essentially, it models a non-linear time-varying system in which every agent's opinion is either a real number or a real-valued vector, and assumes that every agent has a confidence bound defining his/her neighborhood (the set of agents influencing him/her at the given point in time). At every time-step, each agent's belief moves to the arithmetic mean of his/her neighbors' beliefs.

 To cite a few notable results, it was shown in~\cite{dittmer2001consensus} that HK dynamics always converge to a steady state in finite time for every set of initial opinions. Later on, the termination time of the dynamics was studied extensively and it is now known that for a system of $n$ agents having scalar opinions, the \beh{maximum} termination time is at least $\Omega(n^2)$ and at most $O(n^3)$ \cite{bhattacharyya2013convergence},~\cite{mohajer2013convergence},~\cite{wedin2015quadratic}. When the opinions are multidimensional, the best known lower and upper bounds are $\Omega(n^2)$ and $O(n^4)$, respectively \cite{bhattacharyya2013convergence},~\cite{martinsson2016improved}. Other properties of interest such as inter-cluster distance and equilibrium stability were studied in~\cite{blondel20072r} and~\cite{blondel2009krause}. 
 
Even though a number of variants of the HK model have been proposed and analyzed (such as the HK model with stubborn and flexible agents~\cite{fu2014opinion}, the inertial HK model~\cite{chazelle2017inertial}, and continuous-time noisy variants~\cite{wang2017noisy}), very few models, such as the \textit{social HK model}, proposed in \cite{fortunato2005consensus}, the generalized Deffuant-Weisbuch model proposed in \cite{urbig2007communication}, and the \textit{social similarity-based} HK model, proposed in \cite{chen2017opinion}, address an important shortcoming that is central to the original model: the assumption that every agent has \emph{access} to every other agent's opinion (regardless of whether or not he/she is\emph{ influenced} by other agents). 

Such an assumption is questionable, as on large scales a multitude of extrinsic factors such as geographical separation along with differences in culture, nationality, socio-economic background, etc., may drastically reduce the likelihood of two like-minded individuals contacting  each other. To address this issue, the social HK model incorporates a \emph{physical connectivity graph}, denoted by $\Gph$, into the classical HK model. A pair of agents can access each other's opinions if and only if the corresponding vertices are adjacent in $\Gph$.

The social HK model was proposed in \cite{fortunato2005consensus}, which provides a conjecture on the minimum value of the confidence bound required to achieve consensus in the limit \beh{as the number of agents goes to infinity}. Subsequently, \cite{bhattacharyya2015friends} provided an upper bound on the number of time steps in which two agents separated by a minimum distance influence each other. Recently, in~\cite{parasnis2018hegselmann}, we showed that for any incomplete $\Gph$ and any continuous probability density function having the state space as its support, the expected termination time of social HK dynamics is infinity.

This result motivates us to investigate the convergence properties of the social HK model in this paper. We begin by introducing the original HK model, the social HK model, and the associated terminology in Section \ref{sec:formulation}. In Section \ref{sec:terminationtime}, we provide the proof of the aforementioned result on the expected termination time of the dynamics. In Section \ref{sec:convtime}, we 
show that the conditional upper bound on the maximum $\epsilon$-convergence time provided in~\cite{parasnis2018hegselmann} is applicable to a wider class of initial opinion distributions. In Section \ref{sec:arbit_slow_conv}, we 
show that delaying an event that we call \textit{merging} is the only way to indefinitely delay a social HK system's $\epsilon$-convergence to the steady state. We then provide a set of sufficient conditions and another set of necessary conditions for arbitrarily slow merging, and use the necessary conditions to show that the $\epsilon$-convergence time of a complete $r$-partite graph is bounded. We conclude by observing that these conditions are nearly tight under certain assumptions on the initial opinion distribution, and also provide some future directions.

{A subset of the results of this work have also been  reported in our conference  paper~\cite{HKCDC19} (to appear), where we discriminate between consensus and non-consensus states, and provide sufficient conditions for a physical connectivity graph to have an unbounded convergence time in each case.

\beh{\textbf{Notation:}} We denote the set of real numbers by $\R$, the set of positive real numbers by $\R^+$, the set of integers by $\mathbb{Z}$, the set of positive integers by $\mathbb{N}$, and the set $\mathbb N\cup\{0\}$ by $\mathbb{N}_0$. We define $[n]:=\{1,\ldots,n\}$. We use $I$ to denote the identity matrix (of the known dimension). 

We denote the cardinality of a set $S$ by $|S|$, the vector space of column vectors consisting of $n$-tuples of real numbers by $\R^n$, the $\infty$-norm in $\R^n$ by $\norm{\cdot}_\infty$, and the all-one vector and the all-zero vector in $\R^n$ by $\allone_n$ and $\allzero_n$, respectively, dropping the subscripts when the dimension is clear from the context. For a set $S$, $\allone_S$ denotes $\allone_{|S|}$.

\beh{An undirected graph on $n$ vertices is $G=(V,E)$ where  $V$ or $V(G)$ is the set of vertices and $E=E(G)\subseteq V\times V$ is the set of edges, with $(i,j)\in E$ if and only if (iff) $(j,i)\in E$ for $i,j\in V$. If $|V|=n$, we can label the vertices so that $V=[n]$, without loss of generality (w.l.o.g.). For any vector $w\in\mathbb R^n$ and a subset of vertices $V_P\subseteq V$, we let $w_P$ denote the restriction of $w$ to the coordinates specified by $V_P$. Also, for any $l\in [n]$, let $w_{[l]}$ denote the vector $[w_1\,\,\ldots\,\,w_l]^T$. Throughout this work, all the graphs are undirected. We say that $i$ and $j$ are neighbors in $G$, if $(i,j)\in E$ (and hence, $(j,i)\in E$). The set of neighbors of a node $i$ in $G$ is the set $\mathcal N_i:=\{j: (i,j)\in E\}$ and the degree of node $i$ is $d_i:=|\mathcal N_i|$. The adjacency matrix of $G=([n],E)$ is the $n\times n$ binary matrix $A_{adj}$ where $(A_{adj})_{ij}=1$ iff $(i,j)\in E$, and the degree matrix of $G$ is the diagonal matrix $D$ with $D_{ii}=d_i$. We define the normalized adjacency matrix of $G$ to be the matrix $A:=D^{-1}A_{adj}$. The Laplacian of $G$ is defined to be $L:=D-A_{adj}$ and the normalized Laplacian of $G$ is defined to be $N:=D^{-1/2}LD^{-1/2}=I-D^{-1/2}A_{adj}D^{-1/2}$. For two graphs $G_1=([n],E_1)$ and $G_2=([n],E_2)$ on $n$ vertices, we let $G_1\cap G_2=([n],E_1\cap E_2)$. For any subscript $_P$, the normalized adjacency matrix of a graph denoted by $G_P$ is denoted by $A_P$. Finally, a complete graph (or clique) on $n$ vertices is the graph $K_n:=([n],[n]\times [n])$. }

%% file: formulation.tex
\section{PROBLEM FORMULATION}\label{sec:formulation}
\subsection{Original Model}
Consider a network of $n$ agents. For each $k\in \mathbb N$, let $x_i[k]$ be the opinion of the $i$\textsuperscript{th} agent at time $k$. Then the state of the system at time $k$ is defined as $x[k]:= \left[x_1[k]\,\,x_2[k]\,\,\ldots\,\,x_n[k]\right]^T\in\R^n$. Occasionally, we drop the indexing $[k]$ for the state and its associated quantities when the context makes the time index clear. In the original HK model, at time $k$, agents $i$ and $j$ are {neighbors} iff $|x_i[k]-x_j[k]|\leq R$, where $R$, the \emph{confidence bound}, is assumed to be the same for every agent. Thus, the set of neighbors of agent $i$ at time $k$ is: 
$$\mathcal N_i(x[k])=\left\{j\in[n]: |x_i[k]-x_j[k]|\leq R\right\}.$$
Note that $i\in \mathcal N_i$ for all $i\in[n]$. Also, $i$ is a neighbor of $j$ iff $j$ is a neighbor of $i$. Therefore, we can encode all of the information about the influences in the network at time $k$ into an undirected graph, $G_c(x[k])$, which we call the \emph{communication graph} of the network at time $k$. This  $n$-vertex graph has a link between two vertices iff the corresponding agents are neighbors at time $k$. Observe that $G_c(x[k])$ always has a self-loop at each vertex at all times. Finally, at every time instant, every agent's opinion shifts to the average of his/her neighbors' current opinions:
\begin{align}\label{eqn:updaterule}
x_i[k+1] = \frac{\sum_{j\in \mathcal N_i(x[k])}x_j[k]}{|\mathcal N_i(x[k])|}.
\end{align}

This being a \textit{bounded confidence} model, it is possible that an agent does not have any neighbor other than himself/herself, in which case, his/her opinion does not change i.e., $x_i[k+1] =x_i[k]$. Such an agent is said to be \emph{isolated}. 

\subsection{\beh{Modification}}
In the original HK dynamics, if the opinions of any two agents are within a distance of $R$ from each other, then the agents necessarily influence each other. This assumption is relaxed in the social HK model by the introduction of a second graph, as described below.

Let the physical connectivity graph \beh{$G_{ph}=([n],E_{ph})$} be an $n$-vertex graph wherein each vertex represents an agent. Two agents $i$ and $j$ can communicate with each other iff their corresponding vertices are adjacent in $G_{ph}$. Hence, for two individuals to influence each others' opinions, they not only need to be similarly opinionated but also to be physically connected through $G_{ph}$. \beh{Throughout this paper, we assume that $G_{ph}$ is connected, time-invariant, and contains all the self-loops, i.e., $(i,i)\in E_{ph}$ for all $i\in[n]$.}

Observe that in the special case that $G_{ph}$ is a complete graph, no external restrictions are imposed on the interaction between any two agents. This case, therefore, is equivalent to the well-known original model of the last subsection. However, the social HK generalization starts differing from the original model when there is at least one pair of non-adjacent vertices in $G_{ph}$, as will be revealed next.

\subsection{State-Space Representation}
Each of the two models discussed above has the following state-space representation:
\begin{align}\label{eqn:stateevol}
x[k+1] = A\left(x[k]\right) x[k],
\end{align}
where $A(x[k])$ is the normalized adjacency matrix of $G_{ph}\cap G_c(x[k])$. Thus, $\tilde{G}[k]=G_{ph}\cap G_c(x[k])$ is the effective graph or the \emph{influence graph} at time $k$. The original HK model is a special case with $G_{ph}=K_n$, which gives $\tilde G[k]=G_c(x[k])$.

Note the explicit dependence of the state evolution matrix on the state of the system at time $k$. It arises from the dependence of the structure of the communication graph on the agents' opinions at the concerned time instant.

Now, let ${A}_{adj}\left(x[k]\right)$ denote the adjacency matrix of $\tilde G[k]$ and let $\Dt\left(x[k]\right)$ denote its degree matrix. Then
\begin{align*}
x[k+1] = \Dt^{-1}\left(x[k]\right)\cdot \At_{adj}\left(x[k]\right)\cdot x[k]
\end{align*}
which can be expressed more compactly as:
\begin{align}\label{eqn:dynstate}
x[k+1] = \Dt^{-1} \At_{adj} x[k].
\end{align}
In other words, the state evolution matrix is given by $A=D^{-1}A_{adj}$. (We drop the dependencies of these matrices on $x[k]$ for notational simplicity). 

%% file: terminationtime.tex
\section{ANALYSIS OF TERMINATION TIME}\label{sec:terminationtime}
{With the help of standard consensus results such as \cite{lorenz2005stabilization}, one can prove that irrespective of the initial state of the social HK system, its convergence to a steady state is certain, i.e., for every simple physical connectivity graph $\Gph$ and initial state $x[0]=x_0\in\R^n$, the limit $x_\infty(x_0):=\lim_{k\rightarrow\infty}x[k]$ exists. Here, we call $x_\infty(x_0)$ the steady state of the system corresponding to the initial state $x_0$.}

{In this section, using the following definitions, we show that social HK dynamics on an incomplete physical connectivity graph may never attain the steady state in finite time. }
\begin{definition}[Termination Time]\label{def:termtime}
For an initial state $x_0$ and a given physical connectivity graph $\Gph$, the termination time $T(\Gph, x_0)$ is the time taken by the system to reach the steady state corresponding to $x_0$, i.e.:
{$$T(G_{ph},x_0):=\inf\{k\in\mathbb {N}:x[k]=x_\infty(x_0)\}.$$}
\end{definition}
{Next we define the maximum termination time for a given physical connectivity graph.}
\begin{definition}[Maximum Termination Time]\label{def:maxtermtime}
{For a given physical connectivity graph $\Gph$, the \emph{maximum} termination time $T^*(\Gph)$ is the supremum of termination times over all possible initial states:
$$T^*(\Gph):=\sup_{x_0\in\mathbb R^n}T(\Gph,x_0).$$}
\end{definition}

As a special case, it was shown in~\cite{mohajer2013convergence} and~\cite{wedin2015quadratic} that the maximum termination time of the original HK dynamics satisfies {$cn^2 \leq T^*(K_n) \leq C n^3$ asymptotically as $n\rightarrow\infty$ when $d=1$, for some constants $c,C>0$}.

 {We now state a few} propert{ies} of {a class of normalized adjacency matrices that} appear in the state evolution dynamics \eqref{eqn:stateevol}. {These properties} form the basis of our results.

The following lemma is proven in~\cite{martinsson2016improved} as well as \cite{parasnis2018hegselmann}.
\begin{lemma}\label{lemma:diagonalizability} 
	For any undirected graph $\hat{G}$, the normalized adjacency matrix $\hat{A}$ is similar to $I-\hat{N}$ (where $\hat{N}$ is the normalized Laplacian matrix). As a result, $\hat{A}$ is diagonalizable.
\end{lemma}

{The next result provides more information about the spectral properties of the adjacency matrix of a graph, if we have mild additional structures on the graph.}
\begin{lemma}\label{lemma:eigenvalue} 
	Let $\hat{G}$ be an undirected and incomplete graph that is connected and has all the self-loops. Then, if the eigenvalues of the normalized adjacency matrix $\hat A$ (labeled as $\{\lambda_i\}_{i=1}^n$) are ordered such that $|\lambda_1|\geq |\lambda_2|\geq\cdots\geq|\lambda_n|$, we have $1=\lambda_1>|\lambda_2|>0$. Moreover, $\hat A$ has at least one positive eigenvalue besides 1.
\end{lemma}
\begin{proof}
The first part of the result is proven \cite{parasnis2018hegselmann}. Since $\hat A$ is a row-stochastic matrix, we have $\lambda_1=1$. To show that $\hat{A}$ has a positive eigenvalue besides $1$, we have $\sum_{i=1}^{n}\lambda_i=1+\sum_{i=2}^n\lambda_i=\text{Tr}(\hat{A})$, but since $\hat{G}$ is incomplete, $\text{Tr}(\hat{A})>1$. Therefore, $\sum_{i=2}^n\lambda_i>0$ and hence, $\lambda_i>0$ for some $i$.
\end{proof}

We are now ready to show that on average, social HK dynamics on an incomplete graph never terminate.
\begin{proposition}\label{prop:expected}
{Let $x[0]$ be a random vector over $\R^n$ whose distribution induces the Borel-measure $\mu$ on $\R^n$ with $\mu(V)>0$ for any non-empty open set $V\subset \R^n$ (or in other words, the probability density function of $x[0]$ is non-zero almost everywhere). Suppose that $G_{ph}$ is {not a complete graph}}. Then the expected termination time of the dynamics is infinite, i.e., $\mathbb E_{x[0]}[T(\Gph,x[0])]=\infty$.
\end{proposition}
\begin{proof}
		It is sufficient to show that $\Pr(T(\Gph,x[0])=\infty)>0$. Let 
		$S:=\left\{x\in\mathbb R^n:|\max_{i\in[n]} x_i-\min_{j\in[n]} x_j|< R\right\}$.
Note that $S$ is a nonempty open set in $\R^n$ and hence, $\mu(S)>0$. Also, whenever $x[0]\in S$, every agent is within the confidence of every other agent and thus $G_c(x[0])$ is an $n$-clique. Also, from the update rule \eqref{eqn:updaterule}, it {follows} that $\max_{i=1}^n x_i[k]$ is monotonically non-increasing and likewise, $\min_{j=1}^n x_j[k]$ is monotonically non-decreasing {(as functions of time $k$)}. Therefore, the communication graph remains a clique for all $k$, meaning that $\tilde G[k]=G_{ph}$ for all $k\in\mathbb N$. In this case, the dynamics become linear and time-invariant: $A=A(x[k]) = A(x[0])$ and hence, $ x[k]=A^kx[0]$. 

Furthermore, the diagonalizability of $A$ (Lemma \ref{lemma:diagonalizability}) implies that we can write any initial state $x_0\in\R^n$ as a linear combination of the eigenvectors of $A$, i.e., there exist coefficients $c_1(x_0), c_2(x_0), \ldots, c_n(x_0)\in \R$ such that $x_0=c_1(x_0)\allone+\sum_{i=2}^n c_i(x_0)v_i$,
	 where $\allone, v_2,\ldots,v_n$ are  eigenvectors of $A$ corresponding to $\lambda_1=1$, $\lambda_2, \ldots, \lambda_n$. This means that:
$x[k]=c_1(x_0)\allone+\sum_{i=2}^n\lambda^k_ic_i(x_0)v_i$,
for $k\in\mathbb N$. Thus, we have $x_\infty=\lim_{k\rightarrow\infty}x[k]=\alpha\allone$ for some $\alpha\in\R$ because $|\lambda_i|<1$ for $i\not=1$ (by Lemma~\ref{lemma:eigenvalue}).


{Now, consider a random initial vector $x[0]=x_0\sim\mu$. For $k\in\mathbb N$, let us define the event $E_k$ to be the event where the termination time is $k$, i.e., $E_k=\{\omega\mid T(\Gph,x_0(\omega)) = k\}$. Then by our definition of termination time,}
{\begin{align*}
&\Pr(E_{k}\,\mid \,x_0\in S)\cr
&=\Pr\left(\left\{\exists\, x_\infty\in\mathbb R^n: x[\ell]=x_\infty\text{ iff }\ell\geq k\right\}\mid x_0\in S\right)\cr
&\leq\Pr\left(\left\{\exists\, x_\infty\in\R^n: x[k]=x_\infty\right\}\mid \,x_0\in S\right)\cr
&=\Pr\left(\left\{\exists\, \alpha\in\mathbb R: x[k]=\alpha \allone \right\}\mid \,x_0\in S\right)\cr
&=\Pr\left(\left\{\exists\, \alpha\in\mathbb R: (c_1-\alpha)\allone+\sum_{i=2}^n\lambda^k_ic_iv_i=\allzero \right\}\mid \,x_0\in S\right)\cr
&=\Pr\left(\left\{\exists\, \alpha'\in\mathbb R: \alpha'\allone+\sum_{i=2}^n\lambda^k_ic_i(x_0)v_i=\allzero \right\}\mid \,x_0\in S\right).
\end{align*}}
Therefore, {using the fact that the eigenvectors of a diagonalizable matrix are linearly independent and from the fact that $|\lambda_2|>0$ (by Lemma \ref{lemma:eigenvalue}), we have: 
\begin{align}\label{eqn:probbounding}
    &\Pr(E_{k}\,\mid \,x_0\in S)\cr
&\leq \Pr\left(\left\{\lambda_i^kc_i(x_0)=0, i\in[n]\backslash\{1\} \right\}\mid \,x_0\in S\right)\cr
&=\Pr\left(\left\{c_2(x_0)=0, \lambda_i^kc_i(x_0)=0, i\in[n]\backslash\{1,2\} \right\}\mid \,x_0\in S\right)\cr
&\leq\Pr\left(\left\{c_2(x_0)=0\right\}|\,x_0\in S\right).
\end{align}
\noindent Observe ${c_2(x_0)=0}$ only when $x_0\in \spann(\allone,v_3,\ldots,v_n)$ which is a subspace of dimension $n-1$. By the continuity of $x_0$, it follows that $\Pr({c_2(x_0)=0}\mid x_0\in S)=0$.}
%
Using this, and \eqref{eqn:probbounding}, we obtain $\Pr(E_{k}\,\mid \,x[0]\in S)=0$. Therefore, the conditional probability of finite time termination, given that the initial state lies in $S$ is: 
\begin{align*}
\Pr(T(\Gph,x[0])<\infty\mid x[0]\in S)&=\Pr(\cup_{k=0}^\infty E_k \mid x[0]\in S)\\&=0.
\end{align*}
where the last equality follows from the fact that a countable union of infinitely many zero-probability events is also a zero-probability event. We conclude the proof as follows:
\begin{align*}
&\Pr(T(\Gph, x[0])=\infty)\\
&\geq\Pr\left(T(\Gph, x[0])=\infty\,\mid \,x[0]\in S\right)\Pr(x[0]\in S)\\
&=\left(1-\Pr(T(\Gph, x[0])<\infty\,\mid \,x[0]\in S)\right)\Pr(x[0]\in S)\\
&=\Pr(x[0]\in S)>0.\end{align*}\end{proof}
{ In essence, Proposition \ref{prop:expected} states that the expected termination time of the social HK dynamics on any underlying {incomplete $\Gph$} is infinity, which means that there is a continuum of initial states starting from which social HK dynamics never terminate. This shows that the behavior of the HK dynamics over complete graphs is indeed an anomaly.}

%% file: upperbound.tex
\section{BOUNDS ON THE CONVERGENCE TIME}\label{sec:convtime}
{Now that we know that a social HK system may never reach the steady state, the next pertinent question is: {how fast does it approach the steady state?}}

{We begin with a few relevant definitions.}


\begin{definition}[$\epsilon$-Convergence]\label{def:convrate}
{Given a physical connectivity graph $\Gph$, an initial state $x_0$, and $\epsilon>0$, the system is said to have achieved $\epsilon$-convergence at time $N\geq 0$ if its state lies in the $\epsilon$-neighborhood of the steady state corresponding to $x_0$, i.e., $\|x[k]-x_\infty(x_0)\|<\epsilon$,
for all $k\geq N$.}
\end{definition}

{Based on this, we define the $\epsilon$-convergence time as:}
\begin{definition}[$\epsilon$-Convergence Time]\label{def:convtime}
{For a given physical connectivity graph $\Gph$, an initial state $x_0\in\R^n$, and a given $\epsilon>0$, the $\epsilon$-convergence time $k_{\epsilon}(\Gph,x_0)$ is the time taken by the system to achieve $\epsilon$-convergence:}
\begin{align*}
k_\epsilon(\Gph,x_0):&=\inf\left\{N\in\mathbb N:\|x[k]-x_\infty(x_0)\|<\epsilon,\right.\\
&\quad \qquad \left.\text{ for all } k\geq N\right\}.
\end{align*}
\end{definition}
{Similar to $T^*$, we define $k_\epsilon^*$ to be the supremum of $\epsilon$-convergence times for all initial states.}
\begin{definition}[Maximum $\epsilon$-Convergence Time]\label{def:maxconvtime}
For a given physical connectivity graph $\Gph$ and $\epsilon>0$, the maximum $\epsilon$-convergence time $k_{\epsilon}^*(\Gph)$ is the supremum of $\epsilon$-convergence times over all possible initial opinions:
\begin{align}\label{eqn:kepdef}
k_\epsilon^*(\Gph):=\sup_{x_0\in\mathbb R^n}k_\epsilon(\Gph, x_0).
\end{align}

\end{definition}

\subsection{Lower Bound}
{We now provide a lower bound on the maximum $\epsilon$-convergence time $k_\epsilon^*(\Gph)$ in terms of the conductance of $G_{ph}$.} We borrow the definition of conductance from \cite{kannan2004clusterings}.

{Let $G=([n],E)$ be an undirected graph on $n$ vertices. For a subset $S\subset [n]$, let $\partial(S):=\{(i,j)\in E\mid i\in  S, j\in \bar{S}\}$,  where  $\bar{S}=[n]\backslash S$. In words, $\partial S$ represents the set of edges that {connect} $S$ {to} the rest of the graph. Further, let $d(S)$ denote the sum of the degrees of the vertices in $S$. Then we have the following definition.}
\begin{definition}[Conductance] 
The conductance $\phi(G)$ of a graph $G=([n],E)$ is defined as:
{$$\phi(G):=\min_{\substack{S\subset [n]\\S\not=\emptyset}}\frac{|\partial(S)|}{\min\left(d(S),d(\bar S)\right)}.$$}
\end{definition}
The next proposition states that a system whose physical connectivity graph has a low conductance might take a long time to converge to its steady state. See~\cite{parasnis2018hegselmann} for the proof.

\begin{proposition}\label{prop:phi}
{For any incomplete graph $\Gph$ and any given $\epsilon>0$, the maximum $\epsilon$-convergence time of the {social} Hegselmann-Krause dynamics, as defined in \eqref{eqn:kepdef}, satisfies
\begin{align}\label{eqn:lowerbound}
k_\epsilon^*(\Gph)>\frac{\log{\left(\frac{\epsilon{\sqrt 2}}{R}\right)}}{\log\left(1-2\phi(G_{ph})\right)}.\\
\nonumber
\end{align}}
\end{proposition}
\begin{remark}
{Proposition \ref{prop:phi} can be used to compute a lower bound on the maximum $\epsilon$-convergence time in terms of $n$ for graphs whose conductance is known as a function of $n$. For example, the dumbbell graph on $n$ vertices has $\phi=O\left(\frac{1}{n^2}\right)$ \cite{kannan2006blocking} which yields $k_\epsilon^*=\Omega(n^2)$.}
\end{remark}
\subsection{Upper Bound Applicable to a Class of Initial Opinions}

We now show that if the influence graph remains connected and time-invariant until $\epsilon$-convergence to the steady state, then the latter is achieved in $O(n^2\log n\cdot d(G_{ph}))$ {steps, where $d(G_{ph})$ denotes the diameter of $\Gph$}.

\begin{proposition}\label{prop:cond_bound}
Suppose there exist $\epsilon>0$ and an initial state $x_0\in\R^n$ such that the influence graph, $\tilde G[k]$ remains connected and constant in time until $\epsilon$-convergence is achieved. Then with $x_0$ as the initial state, the social HK system achieves $\epsilon$-convergence in $O(n^2\log n\cdot d(\Gph))$ steps.
\end{proposition}
\begin{proof}
Let $\epsilon$ and $x_0$ be as described above, and let $x[0]=x_0$. Then observe that the state evolution until $\epsilon$-convergence can be expressed as $x[k]=A^kx_0$. For $i\in\{2,\ldots,n\}$, let $v_i$ denote an eigenvector of $A$ corresponding to $\lambda_i$. Then, since $A$ is diagonalizable (by Lemma \ref{lemma:diagonalizability}), we have $x_0\in \spann\{\allone, v_2, \ldots, v_n\}$, i.e., $x_0=c_1\allone+ \sum_{i=2}^n c_i v_i$ for some $c_1, c_2, \ldots, c_n\in\mathbb R$. Consequently, $x[k]=c_1\allone+ \sum_{i=2}^n c_i \lambda_i^k v_i$ for $0\leq k\leq k_\epsilon(\Gph,x_0)$.
Therefore, according to Lemma \ref{lemma:diagonalizability}:
\begin{align*}
&\|x[k]-c_1\allone\|=\norm{\sum_{i=2}^n c_i\lambda_i^kv_i}=\norm{A^k\sum_{i=2}^n c_iv_i}\\
                             &=\norm{\Dt^{-1/2}(I_n - \Nt)^k\Dt^{1/2}\sum_{i=2}^n c_iv_i}=\norm{\Dt^{-1/2}M^ky}.
\end{align*}
where $M=I_n-\Nt$ and $y=\sum_{i=2}^n c_i\Dt^{1/2}v_i$. Hence, 
$$\|x[k]-c_1\allone\|\leq\|\Dt^{-1/2}\|_{ind}\cdot\|M^ky\|.$$
Now observe that since $\Dt$ is a positive diagonal matrix,
\begin{align*}
\|\Dt^{-1/2}\|_{ind}&=\max_{i\in [n]}\{(\Dt^{-1/2})_{ii}\}=\left(\min_{i\in [n]} |N_i|\right)^{-\frac{1}{2}}\leq1
\end{align*}
where the inequality is due to the fact that every vertex of $\Gt$ has a self-loop. Thus, 
\begin{equation}\label{eq:M_y}
\|x[k]-c_1\allone\|\leq\|M^ky\|=\sqrt{y^T M^{2k}y},
\end{equation}
because $M$ is symmetric. Next, note that $A=\Dt^{-1/2}M\Dt^{1/2}$ implies that $\Dt^{1/2}\allone, \Dt^{1/2}v_2, \ldots, \Dt^{1/2}v_n$ are the orthogonal eigenvectors of $M$. Therefore $y^T \Dt^{1/2}\allone=0$ and the Courant-Fischer theorem~\cite{golub2012matrix} allows us to bound the expression above as:
\begin{equation}
\sqrt{y^T M^{2k}y}\leq\sqrt{\lambda_2^{2k}y^Ty}=|\lambda_2|^k\|y\|.
\end{equation}
The next step is to bound $\norm{y}$: 
\begin{align}
\norm{y}&=\norm{\Dt^{1/2}\sum_{i=2}^nc_iv_i}\nonumber\\
&\leq  \left(\max_{i\in [n]} \Dt_{ii}\right)^{\frac{1}{2}}\norm{\sum_{i=2}^nc_iv_i}\leq\sqrt n \norm{\sum_{i=2}^nc_iv_i}.
\end{align}
Now, we upper bound 
$\norm{\sum_{i=2}^nc_iv_i}$ as follows:
\begin{align}\label{eq:kappa}
    \norm{\sum_{i=2}^nc_iv_i}&=\norm{x[0]-c_1\allone}\cr
    &\leq\sqrt{n}\max_{j}|x_j[0]-c_1|\cr
    &\stackrel{(a)}\leq\sqrt{n}\left|\max_px_p[0]-\min_qx_q[0]\right|\cr
    &\stackrel{(b)}\leq\sqrt{n}(n-1)R\leq n^{3/2}R.
\end{align}
Here, (a) follows from $\min_q x_q[0]\leq c_1\leq \max_p x_p[0]$, and (b) holds as $\tilde G[0]$ is connected. Combining \eqref{eq:M_y}-\eqref{eq:kappa} yields $\norm{x[k]-c_1\allone}\leq n^2R|\lambda_2|^k$ and using $|\lambda_2|\leq1-\frac{1}{n^2d(G_{ph})}$ (Proposition 2.3 in~\cite{martinsson2016improved}), we get:
$$\norm{x[k]-c_1\allone}\leq n^2R\left(1-\frac{1}{n^2d(G_{ph})}\right)^{k}.$$
Therefore, the condition below ensures $\norm{x[k]-c_1\allone}\leq \epsilon$: 
\begin{align} \label{eq:conv_upper_bound}
    \kappa(\epsilon):=\frac{\log{\frac{\epsilon}{n^2R}}}{\log\left(1-\frac{1}{n^2d(G_{ph})}\right)}\leq k\leq k_\epsilon(\Gph,x_0).
\end{align}
Now, note that if $\tilde G[k]$ remains constant permanently, then $x_\infty(x_0)=c_1\allone$ because $|\lambda_i|<1$ for $i\geq2$ by Lemma \ref{lemma:eigenvalue}. As a result, $k_\epsilon(\Gph,x_0)=\lceil\kappa(\epsilon)\rceil$. On the other hand, if $\tilde G[k]$ varies after $\epsilon$-convergence, then $\max_{i} x_{i}[k]-\min_{j} x_{j}[k]> R$ for $k=k_\epsilon(\Gph,x_0)$. This enforces $k_\epsilon(\Gph,x_0)<\kappa(R/2)$ because otherwise, as per \eqref{eq:conv_upper_bound}, the network would satisfy $\norm{x[k]-c_1\allone}\leq R/2$ and consequently, $\max_i~x_i[k]~-~\min_j~x_j[k]\leq R$ for some $k\leq k_\epsilon(\Gph,x_0)$, which would wrongly imply that $\tilde G[k]=\Gph$ for all $k\in\mathbb N$. Hence, $k_\epsilon(\Gph, x_0)\leq \min(\lceil\kappa(\epsilon)\rceil,\kappa(R/2))$.

Since $\ln\left(1-\frac{1}{n^2d(G_{ph})}\right)\approx -\frac{1}{n^2d(G_{ph})}$ for sufficiently large $n$, we have $\kappa(\epsilon)\approx n^2d(G_{ph})\left(\log\frac{n^2R}{\epsilon} \right)=O(n^2\log n\cdot d(G_{ph}))$. Thus, $k_\epsilon(\Gph,x_0)=O(n^2\log n\cdot d(\Gph))$.
\end{proof}

%% file: arbit_slow_conv_2.tex
\section{Arbitrarily Slow $\epsilon$-Convergence}\label{sec:arbit_slow_conv}
The results in the previous section prompts us to ask: What if the initial state does not enable $\tilde G[k]$ to remain constant in time? In such cases, the convergence time could be unbounded above if the physical connectivity graph has more than three vertices. In other words, it is possible that $k^*_{\epsilon}(G_{ph})=\infty$.

Here is a relevant example from~\cite{bhattacharyya2015friends}. Let $\Gph$ be the path graph on 4 vertices, and let $\mathcal X = \left\{\left[-R,0,R,-(R-\delta)\right]^T\text{ for }\delta\in(0,R/2)\right\}$. Then note that for $x[0]\in\mathcal X$, we have $x_1[1]=-R/2$, $x_2[1] = 0$, $x_3[1] = R/2$ and $x_4[1]=-(R-\delta)$ because at time 1, the sets of neighbours of the first three agents are $\{1,2\}, \{1,2,3\}$ and $\{2,3\}$ respectively. In $\tilde G[1]$, the fourth agent remains disconnected from the first three agents because $R>2\delta$ and the confidence interval of the fourth agent at time 1 is $[\delta-2R, \delta]$. By induction, we can show that $x[k]=\left[-R/2^k, 0,R/2^k,-(R-\delta)\right]^T$ as long as the third and the fourth agents remain outside each others' confidence intervals, i.e., as long as $R/2^k+R-\delta>R$, or equivalently, as long as $k<\log_2 (R/\delta)$. At time $k=\ceil{\log_2 (R/\delta)}$, however, agents 3 and 4 become neighbors. Thus, at $k= \ceil{\log_2 (R/\delta)}$, the influence graph $\tilde G[k]$ is a connected graph satisfying $\max_{i} x_i[k] - \min_{j} x_j[k] = \max\{R/2^k+R-\delta, 2\cdot R/2^{k}\} \leq R$. This implies that $x_\infty(x[0])=c\mathbf 1$ for some $c\in\mathbb R$. Therefore, $\epsilon$-convergence requires $|x_i[k]-c|\leq\epsilon$ for $i\in[n]$. By the triangle inequality, this in turn requires $| x_3[k] - x_4[k] |\leq2\epsilon$ which is not satisfied for $k< \ceil{\log_2 (R/\delta)}$ and $\epsilon<R/2$. Hence, $k_\epsilon(\Gph, x[0]) \geq \ceil{\log_2 (R/\delta)}$. As a result, $k^*_\epsilon(\Gph)\geq\sup_{\delta\in(0,R)} \ceil{\log_2 (R/\delta)}=\infty$.

We can generalize the example above to graphs having more than 4 vertices by choosing the same initial opinions for agents 1 - 4, setting $x_i[0]=x_4[0]$ for $5\leq i\leq n$, and by repeating the above arguments. Therefore, we may state the following lemma without proof.
\begin{lemma}\label{lem:arbit_slow} For every $n\geq 4$, there exists an $n$-vertex physical connectivity graph $\Gph$ such that $k_\epsilon^*(\Gph)=\infty$ for all $\epsilon\in(0,R/2)$.
\end{lemma}

\subsection{Underlying Phenomenon}In the example leading to Lemma \ref{lem:arbit_slow}, $\tilde G[0]$ was a disconnected graph, and we could indefinitely delay the formation of a link between two connected components of this graph so as to make $k_{\epsilon}(\Gph, x[0])$ arbitrarily large. The next proposition will clarify that for any $\Gph$, this is the only way to make $k_{\epsilon}(\Gph, x[0])$ arbitrarily large.

To establish this result, we define two kinds of events that can change the structure of $\tilde G[k]$ during opinion evolution.

\begin{definition}[Link break]
Let $\Gph=(V,\Eph)$. A link break $i--j$ is said to occur at time $k\geq1$ if $i, j\in V$ are such that the nodes $i$ and $j$ are adjacent in $\tilde G[k-1]$ but non-adjacent in $\tilde G[k]$.
\end{definition}
\par Note that a link $(i,j)\in\Eph$ breaks at time $k$ iff $|x_i[k-1]-x_j[k-1]|\leq R$, and $|x_i[k]-x_j[k]|>R$.
\begin{definition}[Merging]
Let $\tilde G[k_0-1]$ be a  disconnected graph for some $k_0\geq 1$, and let $G_1(x[k])=(V_1,E_1(x[k]))$ and $G_2(x[k])=(V_2,E_2(x[k]))$ be two induced subgraphs of $\tilde G[k]$ that are disconnected from each other in $\tilde G[k]$ at time $k_0-1$. Then $G_1$ and $G_2$ are said to merge at time $k_0$ if there exists a pair of agents $(i,j)\in V_1\times V_2$ such that $i$ and $j$ become neighbors at time $k_0$, i.e., $(i,j)\in \tilde E(x[k_0])$.
\end{definition}

Besides merging and link breaks, the only kind of event that can alter the structure of $\tilde G$ is the formation of a link between two agents belonging to the same component of this graph. We call these events \textit{intra-component link formations}.

We now borrow from \cite{bhattacharyya2015friends} the definition of a Lyapunov function called \emph{energy} and that of a related quantity called \emph{active energy}.

\begin{definition}[Energy]
Let $\tilde G[k]=(V[k],E[k])$. The energy of the social HK system at time $k$ is defined as:
$$\mathcal E[k] := \sum_{(i,j)\in E[k]}|x_i[k]-x_j[k]|^2 + \sum_{(i,j)\notin E[k]}R^2.$$
\end{definition}
\begin{definition}[Active energy] Let $\tilde G[k]=(V[k],E[k])$. The active energy of the social HK system at time $k$ is defined as:
$$\mathcal E_{act}[k]: =\sum_{(i,j)\in E[k]}|x_i[k]-x_j[k]|^2.$$
\end{definition}
Note that $0\leq\mathcal E[k]\leq 2{n\choose 2}R^2$ for all $k\in \mathbb N$.
\begin{lemma}\label{link_break} If $i--j$ occurs at time $k+1$ for some $k\in\mathbb N$, then there exist two agents $p,q\in[n]$ such that $p\in \mathcal N_{i}[k]$, $q\in \mathcal N_j[k]$, and $|x_p[k] - x_q[k]|>R$. 
\end{lemma}
\begin{proof}Suppose the lemma is false, i.e., for every pair $(p,q)\in \mathcal N_i[k]\times \mathcal N_j[k]$, we have $|x_p[k] - x_q[k]|\leq R$. Then:
\begin{align*}
|x_i[k+1]&-x_j[k+1]|\\&\leq \max\left\{|\max_{q\in \mathcal N_j[k]}x_q[k]-\min_{p\in \mathcal N_i[k]}x_p[k]|,\right.\\
&\left.\qquad\quad|\max_{p\in \mathcal N_i[k]}x_p[k]-\min_{q\in \mathcal N_j[k]}x_q[k]|\right\}\\
&\leq\max_{(p,q)\in \mathcal N_i[k]\times \mathcal N_j[k]}|x_p[k] - x_q[k]|\leq R.
\end{align*}
The first inequality stems from the fact that HK dynamics are an averaging dynamics and each agent's opinion at any time instant is bounded by the minimum and the maximum of his/her neighbors' opinions at the previous time instant.
The last inequality above implies that agents $i$ and $j$ are neighbors at time $k+1$, thereby contradicting the fact that the link $(i,j)$ breaks at time $k+1$.
\end{proof} 

Next, we need to establish that only finitely many link breaks can occur in any opinion evolution process.
\begin{lemma}\label{number_breaks} 
    The total number of link breaks during the entire process of opinion evolution is $O(n^5)$ regardless of the structure of $\Gph$ and the initial state $x[0]\in\mathbb R^n$. 
\end{lemma}
\begin{proof} 
Based on Proposition 1 of \cite{bhattacharyya2015friends}, we have:
\begin{equation} \label{eq:energy_decrement}
    \mathcal E[k] - \mathcal E[k+1] \geq (1-|\lambda_k|^2) \mathcal E_{act}[k]
\end{equation}
for $k\in \mathbb N$, where 
\[\lambda_k:=\{\max|\lambda|: \lambda\not=1 \mbox{ is an eigenvalue of $A[k]$}\},\]
and if we let $d_{\mbox{eff}}(G)$ be the largest diameter of any connected component of the graph $G$, we have the lower bound 
\begin{equation}\label{mart}
    1-|\lambda_k|^2\geq \frac{3}{2n^2d_{\mbox{eff}}(\tilde{G}[k])}\geq \frac{3}{2n^3},
\end{equation} 
which was derived in \cite{martinsson2016improved}. Here, we derive a lower bound on the active energy. Let $i,j\in [n]$ and suppose $i--j$ occurs at time $k+1$ for some $k\in\mathbb N$. Then by Lemma \ref{link_break}, we can find two agents $p,q\in[n]$ such that $p\in \mathcal N_{i}[k]$, $q\in \mathcal N_j[k]$, and $|x_p[k] - x_q[k]|>R$. Therefore, by the definition of active energy, we have
\begin{align}\label{active_bound}
&\mathcal E_{act}[k] \\
&\geq |x_p[k]-x_i[k]|^2+|x_i[k]-x_j[k]|^2+|x_j[k]-x_q[k]|^2\nonumber\\
&\geq \frac{1}{3} \left(|x_p[k]-x_i[k]|+|x_i[k]-x_j[k]|+|x_j[k]-x_q[k]|\right)^2\nonumber\\
&\geq\frac{1}{3}|x_p[k]-x_q[k]|^2> R^2/3\nonumber,
\end{align}
 where the second and the third inequalities follow from the Cauchy-Schwarz and the triangle inequalities, respectively. 

Combining \eqref{eq:energy_decrement}, \eqref{mart} and \eqref{active_bound} yields:
\begin{equation}\label{energy_decrement}
\mathcal E[k] - \mathcal E[k+1] \geq \frac{R^2}{2n^3}
\end{equation}
which means that the energy of the network decreases every time a link breaks and the decrement corresponding to each link break is at least ${R^2}/{2n^3}$. Since $\mathcal E[0]\leq n^2R^2$ and $\mathcal E[k]\geq 0$ for $k\in\mathbb N$, the maximum possible number of link breaks that can ever occur is at most
$\frac{n^2R^2}{{R^2}/{2n^3}}=O(n^5)$.

\end{proof}
The next lemma bounds the maximum possible time interval between two consecutive link breaks under the condition that no new link is formed during this interval. 
\begin{lemma}\label{interval} Let $G_P=(V_P, E_P)$ be a connected component of $\tilde G[k_0]$ at some time $k_0\geq 0$. Suppose (i) no link break occurs between any two agents of $G_P$ until time $k_1>k_0$, (ii) one or more link breaks occur within $G_P$ at time $k_1$, and (iii) no new edge is formed  between a node in $V_P$ and another node (in $[n]$) during the time interval $(k_0, k_1)$. Then $k_1 - k_0 = O(n^3 \log n)$.
\end{lemma}
\begin{proof}
First, observe that for every $k$ in the range $k_0\leq k \leq k_1$, there exist two agents $p_k, q_k\in V_P$ such that $|x_{p_k}[k] - x_{q_k}[k]|>R$. If this were false for some $k'\in \{k_0,k_0+1,\ldots,k_1\}$, we would have $\max_{i\in V_P}x_i[k'] - \min_{j\in V_P}x_j[k']\leq R$. Since the difference $D[k]:=\max_{i\in V_P}x_i[k]-\min_{j\in V_P}x_j[k]$ is monotonically non-increasing in $k$, this would imply that every agent of $G_P$ remains within the confidence of every other agent of $G_P$ for all $k\geq k'$, thus contradicting the occurrence of link breaks within $G_P$ at time $k_1$. Hence, $D_e:=D[k_1-1]>R$. 

Next, by assumptions (i) and (iii), the constant graph $G_P$ remains a connected component of $\tilde G[k]$ during the interval $(k_0, k_1)$. Let $c\mathbf 1_{V_P}$ denote the steady state that we would have associated with the original network if $V_P$ were its vertex set and $x_P[k_0]$ (where $x_P$ denotes the restriction of $x$ to the coordinates specified by $V_P$) were its initial state, $y[0]$. In this hypothetical scenario, $\tilde G[k]$ would achieve $D_e$-convergence by time $\Delta:=k_1-1-k_0$ because $\min_iy_i[\Delta]\leq c\leq \max_{i}y_i[\Delta]$ would yield: 
\begin{align*}
\max_i |y_i[\Delta]-c|&= \max\left(c-\min_{i}y_i[\Delta],\max_{i}y_i[\Delta]-c\right)\cr 
&\leq \max_iy_i[\Delta]-\min_iy_i[\Delta]\cr 
&=\max_{i}{x_P}_i[\Delta+k_0]-\min_{j}{x_P}_j[\Delta+k_0]\cr 
&=D_e.
\end{align*}

Moreover, we would have $\tilde G[k]=G_P$ for $k\in (0, \Delta]$. Therefore, for $\epsilon>0$, any $\epsilon$-convergence that would occur by time $\Delta$, would occur in $O(|V_P|^3\log |V_P|)=O(n^3\log n)$ steps. Since $D_e>R>0$, the last paragraph implies that $\Delta=O(n^3\log n)$, thus completing the proof.
\end{proof}

We are now ready to show that merging is unavoidable if we desire arbitrarily slow $\epsilon$-convergence to the steady state.
\begin{proposition}\label{merging_must} In social HK dynamics, all the link breaks and intra-component link formations always occur in $O(n^8\log n)$ time steps. Hence, if there exists $\epsilon>0$ such that $k^*_\epsilon(\Gph)=\infty$, then there exists a set $\mathcal X_0\subset R^n$ such that whenever $x[0]\in\mathcal X_0$, merging occurs at least once during the process of opinion evolution. 
\end{proposition}
\begin{proof}
Let $\epsilon>0$ be such that $k_\epsilon^*(\Gph)=\infty$. Consider an arbitrary initial state $x[0]\in \R^n$. Consider the following two cases in the evolution of the dynamics: Case 1: no link formation ever takes places. Then by Lemma \ref{number_breaks}, we know that at most $O(n^5)$ links break in the opinion evolution process, and by Lemma \ref{interval}, the maximum possible time interval between two consecutive link breaks is $O(n^3\log n)$. Therefore, the time at which the last link breaks is at most $O(n^8\log n)$. After this point in time, the structure of $\tilde G[k]$ never changes. Therefore, for any $\epsilon>0$, it takes $O(n^2\log(n)d(\Gph))$ additional time steps to achieve $\epsilon$-convergence. Hence, $k_{\epsilon}(\Gph,x[0])=O(n^8\log n)+O(n^2\log(n)d(\Gph))=O(n^8\log n)$.

Now, consider Case 2: at least one new link is formed during the evolution of the dynamics (from the initial state $x[0]$) but no merging ever occurs. For $r\in\mathbb N\backslash\{0\}$, let $t_r$ denote the time at which the $r$-th set of simultaneous link breaks occur and w.l.o.g., suppose the first link formation occurs at a time $k'\in\{t_l, t_l +1, \ldots, t_{l+1}\}$ for some $l\in\mathbb N\backslash\{0\}$. Let $(i,j)$ denote this new link. Since no merging occurs, $(i,j)$ is formed \emph{within} some connected component $G'$ of $\tilde G$. Thus, we have $|x_i[k'-1]-x_j[k'-1]|>R$ and $|x_i[k']-x_j[k']|\leq R$. Also, no link formation or link break during the time interval $[t_l, k'-1]$ implies that $G'$ is a connected component of $\tilde G[k]$ for all $k\in[t_l,k'-1]$. In other words, the influence graph has a connected component that remains constant during the time interval $[t_l, k'-1]$. Therefore, the arguments used in the proof of Proposition \ref{prop:cond_bound} can be repeated to show that $k'-t_l=O(n^2\log n\cdot d(G_1))=O(n^3\log n)$. 

We can repeat the arguments used in the preceding paragraph for subsequent link formations. 

Next, we estimate the maximum number of link formations that can occur in any opinion evolution process. Note that there are at most $n^2$ links in an $n$-vertex graph. So, it may appear that at most $n^2$ link formations can occur. However, every link break gives rise to the possibility of a link formation. Therefore, the maximum number of link formations is $O(n^5)+n^2=O(n^5)$. Hence, if all the link breaks and intra-component link formations were to occur one after the other, then by Lemma \ref{number_breaks}, all of these events would occur in $O(n^5\cdot n^3\log n+ n^5\cdot n^3\log n)=O(n^{8}\log n)$ steps. On the other hand, it is also possible that some of these events occur simultaneously, so that the last of them occurs even sooner. After all the link breaks and link formations, however, the structure of $\tilde G$ remains constant and $\epsilon$-convergence is achieved in $O(n^3\log n)$ additional steps. Thus,  $k_{\epsilon}(\Gph,x[0])=O(n^{8}(\log n))$ in Case 2. 

Finally, since $k^*_\epsilon(\Gph)=\infty$, there exists a set of initial states $\mathcal X_0\subset \R^n$ that do not belong to the above two cases, i.e., a merging event occurs during the evolution of the dynamics started at those initial states. 
\end{proof}
\subsection{Sufficient Conditions for Arbitrarily Slow Merging}
Since the results of the previous subsection imply that arbitrarily slow merging between two components of $\Gph$ is necessary as well as sufficient for arbitrarily slow $\epsilon$-convergence, it is essential to analyze the concept of arbitrarily slow merging in order to better understand the latter concept. For this purpose, we provide conditions on the components of $\Gph$ that ensure that the time at which the corresponding components of $\tilde G[k]$ merge is unbounded when $\epsilon$ is sufficiently small. 
\begin{figure}[h]
 \centering
 \includegraphics[scale=0.5]{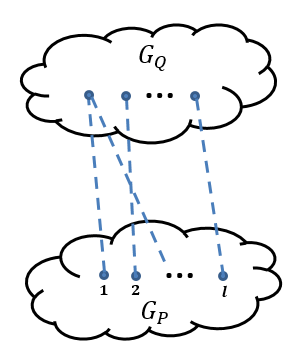}
 \caption{Illustration for Proposition \ref{prop:positive_evector}}
\end{figure}
\begin{proposition}\label{prop:positive_evector}
Suppose $G_{P_0}$ and $G_{Q_0}$ are two subgraphs of $\Gph$ induced by the disjoint vertex sets $V_P\subset V$ and $V_Q\subset V$, respectively. Also, suppose the following conditions hold:
\begin{enumerate}
\item $G_{P_0}$ is a connected graph.
\item W.l.o.g., let $1,2,\ldots,l\in V_P$ be the nodes of $G_{P_0}$ that are adjacent to $V_Q$. Then $A_{P_0}$ has an eigenvalue $\lambda$ such that $0<\lambda<1$, and there exists a corresponding eigenvector $v$ such that $v_1,v_2,\ldots,v_l$ are all of the same sign, i.e., $v_iv_j>0$ for all $i,j\in[l]$.
\end{enumerate}
Let $G_P=G_P[k]$ and $G_Q=G_Q[k]$ denote the subgraphs of $\tilde G[k]$ induced by $V_P$ and $V_Q$, respectively, and for $x_0\in\R^n$, let $k_{M}(x_0,V_P,V_Q)$ denote the time at which $G_{P}$ and $G_{Q}$, merge for the first time, under the condition that the initial state is $x_0$. Then there exists a set $\mathcal X_M$ of initial states such that $k_M(x_0,V_P,V_Q)<\infty$ for all $x_0\in\mathcal X_M$ (i.e., merging occurs) and $\sup_{x_0\in\mathcal X_M} k_M(x_0, V_P, V_Q)=\infty$.
\end{proposition}
\begin{proof}W.l.o.g., let $V_P=\{1,2,\ldots, p\}$ 
and $V_Q= \{p+1,p+2,\ldots, p+q\}$ for some $p,q\in[n]$. Scale $v$ (outlined in condition (2) of the proposition) appropriately so as to satisfy $v_i<0$ for $1\leq i\leq l$ and $\max_i v_i - \min_j v_j \leq R$. Let $v_0:=-\max_{i=1}^l v_i$. Consider $\mathcal X_M=\{z\in\mathbb R^n: z_P=v, z_Q = (R-\delta)\mathbf 1_{q} \text{ for some }\delta\in(0,v_0)\}$. Then observe that if $x[0]\in\mathcal X_M$, the range of allowed values of $\delta$ ensures that $\tilde G[0]$ is a disjoint union of $G_{P_0}$, $G_{Q_0}$, and possibly some other connected components. This is because all the potential neighbours of $V_Q$ in $V_P$, namely the nodes $1,2,\ldots l$, are outside the confidence interval $[-\delta, 2R-\delta]$ of every agent in $V_Q$, and because $\max_i x_{Pi}[0] - \min_j x_{Pj}[0] \leq R$ implies that $G_{P_0}$ is an induced subgraph of $\tilde G[0]$. Also, note that $\max_i x_{Pi}[0] - \min_j x_{Pj}[0] \leq R$ enforces $x_P[1]=A_{P_0} x[0] = \lambda v$. 

Now, $\lambda>0$ implies that $x_{Pi}[1]=\lambda v_i<0$ for $i\in [l]$. Therefore, $G_P$ and $G_Q$ are also disconnected from each other in $\tilde G[1]$ provided $\delta<\lambda v_0$. Similarly, for $k\geq 1$, we have $x_{P}[k]=\lambda^k v$ implying that $G_P$ and $G_Q$ remain disconnected from each other as long as $\delta<\lambda^k v_0$, i.e., for $k<\log_{1/\lambda}(v_0/\delta)$. However, since $\lambda<1$, a time is reached when $k=\ceil{\log_{1/\lambda}(v_0/\delta)}$ and consequently, $k_M(x[0],V_P,V_Q)=\ceil{\log_{1/\lambda}(v_0/\delta)}<\infty$ because the agent having the opinion $\max_{i=1}^l x_{Pi}[k]=-\lambda^k v_0$ enters the confidence interval $[-\delta, 2R-\delta]$ of its potential neighbor(s) in $G_Q$. Therefore, $\sup_{x[0]\in\mathcal X_M}k_M(x[0],V_P,V_Q)=\sup_{\delta\in(0,v_0)}\ceil{\log_{1/\lambda}(v_0/\delta)}=\infty$.
\end{proof}




\subsection{Necessary Conditions for Arbitrarily Slow Merging}


Having seen a set of sufficient conditions for arbitrarily slow merging, we now move on to present the necessary conditions for a pair of subgraphs of the influence graph to exhibit this property. To be precise, we ask:
given two subgraphs $G_P[k]$ and $G_Q[k]$ of the influence graph $\tilde G[k]$, can we find a set of necessary conditions for $G_P$ and $G_Q$ to be able to merge at an arbitrary time $\kappa \in\mathbb N$? The main result of this subsection (Proposition $\ref{prop:final}$) answers this question. 

We begin with a few technical lemmas.

\begin{lemma}\label{lem:gamma_lemma}
Consider a vector subspace $U\subset \R^n$ such that for every $v\in U\setminus\{0\}$, we have $v_iv_j<0$ for some $i,j\in [l]$. Further, define $\phi:\R^n\setminus\{0\}\rightarrow\R$ as 
$$\phi(v)=\min\left(\left|\frac{\min_{i\in[l]}v_i}{\max_{i\in[l]}v_i}\right|,\left|\frac{\max_{i\in[l]} v_i}{\min_{i\in[l]} v_i}\right|\right).$$
 Then there exists a constant $\gamma>0$ such that $\phi(v)\geq\gamma$ for all $v\in U\setminus\{0\}$. 
\end{lemma}
\begin{proof}
Since $\phi(\lambda v)=\phi(v)$ for all $\lambda\in\R\setminus\{0\}$ and $v\in\R^n\setminus\{0\}$, it suffices to prove the lemma for $v\in D:=U\cap\{v\in\R^n:\norm{v}=1\}$.

Observe that $D$ is the intersection of the unit ball and a vector subspace of $\R^n$. Hence, it is a compact set. Next, since $\phi$ is continuous on $D$, we know that $\inf_{v\in D}\phi(v)$ is attained because $D$ is a compact set, i.e., $U^*:=\argmin_{v\in D}\phi(v)$ exists and is well defined. Hence, for all $v\in D$ and any $u\in U^*$, we have
$\phi(v)\geq \min_{v\in D}\phi(v)=\phi(u)>0,$ which follows from the assumption of the lemma enforcing $\max_{i\in [l]}u_i> 0> \min_{i\in [l]}u_i $.
\end{proof}

\begin{lemma}\label{sign_lemma}
Let $v\in\mathbb R^l$ satisfy $v_iv_j<0$ for some $i,j\in[n]$, 
 and let $\gamma\in\mathbb R$ be any constant satisfying $0<\gamma\leq|\max_iv_i|/|\min_iv_i|$.
Then for any $u\in\mathbb R^l$, either
\begin{equation}\label{first_statement}
\max_{i}(v+u)_i>0\text{ and }\left|\frac{\max_{i}(v+u)_i}{\min_{i}(v+u)_i}\right|\geq\gamma'
\end{equation}
or
\begin{equation}\label{second_statement} 
\max_{i}(v-u)_i>0\text{ and }\left|\frac{\max_{i}(v-u)_i}{\min_{i}(v-u)_i}\right|\geq\gamma'
\end{equation}
where $\gamma':=\frac{\gamma}{\gamma+2}$.
Moreover, if \eqref{first_statement} holds and $\min_i(v-u)_i<0$, then
\begin{equation}\label{alpha_zeta}
{\max_{i}(v+u)_i}\geq\frac{\gamma|{\min_{i}(v-u)_i}|-\max(0,\max_{i}(v-u)_i)}{\gamma+1}.
\end{equation} 
$$ $$
\end{lemma}
\begin{proof}
We first prove that either \eqref{first_statement} or \eqref{second_statement} holds. Before we begin, observe that $\gamma>0$ implies $\gamma'<\min(1,\gamma)$.

Now, suppose neither \eqref{first_statement} nor \eqref{second_statement} is true. However, we know that $\max_{i}(v+u)_i+\max_{i}(v-u)_i\geq\max_i[(v+u) + (v-u)]_i=2\max_i v_i>0$. As a result, either $\max_i(v+u)_i>0$ (in which case $|\max_i(v+u)_i|<\gamma'|\min_i(v+u)_i|$), or $\max_i(v-u)_i>0$ (in which case $|\max_i(v-u)_i|<\gamma'|\min_i(v-u)_i|$). By implication, there exists a constant $\epsilon\in (0,\gamma')$ such that $P_1\leq\epsilon M_1$ and $P_2\leq\epsilon M_2$, where we define $P_1:=\max(0,\max_{i}(v+u)_i)$, $P_2:=\max(0,\max_{i}(v-u)_i)$, $M_1:=|\min_{i}(v+u)_i|$ and $M_2:=|\min_{i}(v-u)_i|$.

Now, three cases arise.

\emph{Case 1}: $(\min_i(v+u)_i)(\min_i(v-u)_i)\neq0$ and either $\min_i(v+u)_i>0$ or $\min_i(v-u)_i>0$. Suppose $\min_i(v+u)_i=M_1>0$. Then $\max_i(v+u)_i/\min_i(v+u)_i\geq1\geq\gamma'$, thus contradicting the inequality $P_1\leq\epsilon M_1$ and thereby proving the first part of the lemma. The subcase $\min_i(v-u)_i>0$ is handled similarly.

\emph{Case 2}: Either $\min_i(v+u)_i=0$ or $\min_i(v-u)_i=0$. Suppose $\min_i(v+u)_i=M_1=0$. If $\max_i(v+u)_i>0$, then we have $\epsilon M_1=0<P_1$, which again results in a contradiction and establishes the first part of the lemma. On the other hand, if $\max_i(v+u)_i=0$, then it follows that $u=-v$. Consequently, the assumptions made by the lemma lead to the following: $\max_i(v-u)_i=2\max_i v_i>0$, and $\max_i(v-u)_i\geq 2\gamma|\min_i v_i|=\gamma|\min_i (v-u)_i|>\gamma'|\min_i (v-u)_i|$. These inequalities establish \eqref{second_statement} and hence prove the first assertion of the lemma. The subcase $\min_i(v-u)_i=0$ is handled similarly.

\emph{Case 3}: $\min_i(v+u)_i=-M_1<0$ and $\min_i(v+u)_i=-M_2<0$. Observe that 
\begin{align}\label{p1p2}
2\max_iv_i&=\max_i\{(v+u) + (v-u)\}_i\nonumber\\
&\leq \max_{i}(v+u)_i+\max_{i}(v-u)_i\nonumber\\
&\stackrel{(a)}\leq P_1 + P_2\leq \epsilon(M_1+M_2).
\end{align}
Also,
\begin{align}\label{p1m2} 
2\min_i v_i&= \min_i\{(v+u) + (v-u)\}_i\nonumber\\
&\leq\max_{i}(v+u)_i+\min_{i}(v-u)_i\nonumber\\
&\stackrel{(b)}\leq P_1-M_2\leq \epsilon M_1-M_2. 
\end{align}
Similarly,
\begin{equation}\label{p2m1}
2\min_i v_i\leq\epsilon M_2-M_1.
\end{equation}

Within this case, two subcases arise. Subcase 1: Suppose both $\epsilon M_1-M_2<0$ and $\epsilon M_2-M_1<0$. In other words, $\epsilon<\eta<1/\epsilon$, where we define $\eta:=M_1/M_2$. Consider the inequality $\epsilon M_1-M_2<0$ first. Along with \eqref{p1m2}, it implies: $|\min_i v_i|\geq0.5|M_2-\epsilon M_1|$. Likewise, \eqref{p1p2} and the assumption  $\max_i v_i>0$ imply: $|\max_i v_i|\leq 0.5\epsilon(M_1+M_2)$. Combining these inequalities with the assumption $\max_i v_i\geq\gamma|\min_i v_i|$ yields:
\begin{align}\label{ineq_one}
\frac{\epsilon(M_1+M_2)}{M_2-\epsilon M_1}\geq\gamma.
\end{align}
Similarly, the subcase inequality $\epsilon M_2-M_1<0$, leads to:
\begin{align}\label{ineq_two}
\frac{\epsilon(M_1+M_2)}{M_1-\epsilon M_2}\geq\gamma.
\end{align}
We express \eqref{ineq_one} and \eqref{ineq_two} in terms of $\eta$ as \eqref{last}-(a) and \eqref{last}-(b) respectively:
\begin{equation}\label{last}
\frac{\gamma-\epsilon}{\epsilon(1+\gamma)}\stackrel{(a)}\leq\eta\stackrel{(b)}\leq\frac{\epsilon(1+\gamma)}{\gamma-\epsilon},
\end{equation}
which is possible only if $\gamma-\epsilon\leq\epsilon(1+\gamma)$, i.e., only if $\epsilon\geq\frac{\gamma}{\gamma+2}=\gamma'.$ This contradicts that $\epsilon\in(0,\gamma')$, thus establishing the first assertion of the lemma.

Finally, we have Subcase 2: $\epsilon M_2-M_1\geq0$ or $\epsilon M_1-M_2\geq0$. We assume the former w.l.o.g. Then $\eta\leq\epsilon< \gamma'<1<1/\epsilon$. Hence $\epsilon M_1-M_2<0$, implying \eqref{last}-(a) again. On eliminating $\eta$ by using the observation $\eta\leq\epsilon$, we obtain $(\gamma+1)\epsilon^2+\epsilon-\gamma\geq 0$. Since $\epsilon$ is positive by assumption, this inequality requires $\epsilon\geq \frac{\gamma}{\gamma+1}\geq\gamma'$ which contradicts $\epsilon\in(0,\gamma')$, thereby proving that either \eqref{first_statement} or \eqref{second_statement} holds. 

For the second part, given that \eqref{first_statement} holds, we have $\max_i(v+u)_i=P_1>0$. Note that if $P_1\geq M_2=|\min_i(v-u)_i|$, then \eqref{alpha_zeta} follows from $\gamma>0$. So, suppose that $P_1<M_2$. Then \eqref{p1m2}-(b) implies that $2|\min_i v_i|\geq M_2-P_1$. Likewise, $\max_i v_i>0$ and \eqref{p1p2}-(a) together imply that $2|\max_i v_i|\leq P_1+P_2$. Combining these inequalities with the lemma assumption $\max_iv_i\geq\gamma|\min_iv_i|$ yields $
\frac{P_1+P_2}{M_2-P_1}\geq\gamma$, rearranging which we obtain:
$$P_1\geq\frac{\gamma M_2-P_2}{\gamma+1}$$
which is equivalent to \eqref{alpha_zeta}.
\end{proof}

For the next result, we will need to consider a normalized adjacency matrix, suitably combine its eigenvectors associated with repeated eigenvalues, and carefully account for the sign flips arising from powers of negative eigenvalues. 
For this purpose, we introduce the \textit{Elimination Method}.

\subsection*{The {Elimination Method}}
Let $\{\lambda_i\}_{i=1}^{\tau}\subset\R$, $\{u_{\boldsymbol i}\}_{i=1}^\tau\subset\R^n$ and $l\in[n]$ be fixed, and let
$
    S_{j}[k]:=\sum_{i=1}^{\tau}\lambda_i^k u_{\boldsymbol ij}
$
for all $k\in\mathbb N$ and $j\in [l]$. Then the Elimination Method is as follows:
\begin{enumerate}
    \item Find a minimal set of real numbers, $\{\mu_i\}_{i=1}^m$, such that $\{\lambda_i\}_{i=1}^{n}=\{\pm\mu_i\}_{i=1}^m\cup\{0\}$, and $\mu_1>\cdots>\mu_m>0$.
    \item For each $i\in [m]$, find $\nu_i, \sigma_i\in [\tau]$ satisfying $\lambda_{\nu_i}=\mu_i=-\lambda_{\sigma_i}$, and define $u_{\boldsymbol i}^+:=u_{\boldsymbol {\nu}_i}$ and $u_{\boldsymbol i}^-:=u_{\boldsymbol {\sigma}_i}$. If no such $\nu_i$ (respectively, $\sigma_i$) exists, then set $u_{\boldsymbol i}^+=0$ (respectively, $u_{\boldsymbol i}^-=0)$.
    \item For $i\in [m]$, define $\alpha_i=\max_{j=1}^l|u_{\boldsymbol ij}^{+}+u_{\boldsymbol ij}^{-}|$ and $\zeta_i=\max_{j=1}^l|u_{\boldsymbol ij}^{+}-u_{\boldsymbol ij}^{-}|$. Further, define $v_{\boldsymbol i}$ by: $v_{\boldsymbol i}:=\left(u_{\boldsymbol i}^++u_{\boldsymbol i}^-\right)/\alpha_i$ if $\alpha_i\neq0$ and $v_{\boldsymbol i}:=0$ otherwise. Likewise, let $z_{\boldsymbol i}:=\left(u_{\boldsymbol i}^+-u_{\boldsymbol i}^-\right)/\zeta_i$ if $\zeta_i\neq0$ and let $z_{\boldsymbol i}:=0$ otherwise.
\item If $\alpha_i=\zeta_i=0$, discard $\mu_i$ from $\{\mu_j\}_{j=1}^m$, decrement the value of $m$ by 1, and re-enumerate $\{\mu_j\}_{j=1}^m$ so that $\mu_1>\mu_2>\cdots>\mu_m$.
\end{enumerate}
As a result of this procedure, we have the following relations for all $j\in[l]$:
\begin{align}
    S_j[k] &=\sum_{i=1}^m\alpha_i\mu_i^kv_{\boldsymbol ij}\quad \forall\,\,k\in\mathbb N:k\text{ is even},\\
    S_j[k] &=\sum_{i=1}^m\zeta_i\mu_i^kz_{\boldsymbol ij}\quad \forall\,\,k\in\mathbb N:k\text{ is odd.}
\end{align}
The vectors, $\{v_{\boldsymbol{i}}\}_{i=1}^m$ and $\{z_{\boldsymbol{i}}\}_{i=1}^m$ will be called even-$k$ vectors and odd-$k$ vectors respectively.

 We now recast Lemma \ref{sign_lemma} into a more useful form.
\begin{lemma}\label{lem:sign_extended}
    Let $\{\lambda_i\}_{i=1}^{\tau}\subset\R$, $\{u_{\boldsymbol i}\}_{i=1}^\tau\subset\R^n$ and $l\in[n]$ be fixed, and let
$
    S_{j}[k]:=\sum_{i=1}^{\tau}\lambda_i^k u_{\boldsymbol ij}
$. Suppose that
    for every pair $(\lambda_i,u_{\boldsymbol i})$ satisfying $0<\lambda_i<1$ and $(u_{\boldsymbol i})_{[l]}\neq 0$, we have $\max_{p\in [l]}u_{\boldsymbol ip}>0$, $\min_{p\in [l]}u_{\boldsymbol ip}<0$ and $|\max_{p\in [l]}u_{\boldsymbol ip}|/|\min_{q\in[l]}u_{\boldsymbol iq}|\geq\gamma_0$, where $\gamma_0\in\R^+$ is a constant.
    Further, for each $i\in [m]$, let $p(i):=\argmax_{j=1}^l v_{\boldsymbol ij}$ and $\tilde p(i):=\argmax_{j=1}^l z_{\boldsymbol ij}$, where $v_{\boldsymbol i}$ and $z_{\boldsymbol i}$ are given by the Elimination Method. Then for every $i\in[m]$, we have $$\max\left(v_{\boldsymbol ip(i)},z_{\boldsymbol i\tilde p(i)}\right)\geq\hat\gamma_0,$$
    where $\hat\gamma_0:=\frac{\gamma_0}{2+\gamma_0}.$ Furthermore, if $v_{\boldsymbol ip(i)}<\hat\gamma_0$ (respectively, $z_{\boldsymbol i\tilde p(i)}<\hat\gamma_0$), then $\zeta_i/\alpha_i\geq\hat\gamma_0$ (respectively $\alpha_i/\zeta_i\geq\hat\gamma_0$).
\end{lemma}

\begin{proof}
Consider any $i\in [m]$ and let $q(i):=\argmin_{j=1}^l v_{\boldsymbol ij}$ and $\tilde q(i):=\argmin_{j=1}^l z_{\boldsymbol ij}$.

 Now, two possibilities arise: either $(u_{\boldsymbol i}^+)_{[l]}= 0$ or $(u_{\boldsymbol i}^+)_{[l]}\neq 0$. If $(u_{\boldsymbol i}^+)_{[l]}=0$, then $(v_{\boldsymbol i})_{[l]} = -(z_{\boldsymbol i})_{[l]}$. Hence, either $v_{\boldsymbol ip(i)}\geq|v_{\boldsymbol iq(i)}|$ or $z_{\boldsymbol i\tilde p(i)}\geq|z_{\boldsymbol i\tilde q(i)}|$. Since $\max\left(\max_{f\in [l]}|v_{\boldsymbol if}|,\max_{f\in [l]}|z_{\boldsymbol if}|\right)=1\geq\hat\gamma_0$ due to the Elimination Method, we have $\max\left(v_{\boldsymbol ip(i)},z_{\boldsymbol i\tilde p(i)}\right)\geq\hat\gamma_0$. 

On the other hand, if $(u_{\boldsymbol i}^+)_{[l]}\neq0$, then $\mu_i=\lambda_{\nu_i}$ for some $\nu_i\in[\tau]$. Hence, $\max_{f\in[l]}u_{\boldsymbol if}^+>0$ and $|\max_{f\in[l]}u_{\boldsymbol if}^+|/|\min_{f\in[l]}u_{\boldsymbol if}^+|\geq\gamma_0$. In the light of Lemma \ref{sign_lemma}, this implies that either
\begin{equation}\label{first_thing}
\alpha_iv_{\boldsymbol ip(i)}>0\text{ and }|v_{\boldsymbol ip(i)}|\geq\hat\gamma_0|v_{\boldsymbol iq(i)}|,
\end{equation}
or
\begin{equation}\label{second_thing} 
\zeta_iz_{\boldsymbol i\tilde p(i)}>0\text{ and }|z_{\boldsymbol i\tilde p(i)}|\geq\hat\gamma_0|z_{\boldsymbol i\tilde q(i)}|.
\end{equation}
If \eqref{first_thing} holds, then as a result of the Elimination Method, $1=\max_{f\in[l]}|v_{\boldsymbol if}|=\max\left(|v_{\boldsymbol ip(i)}|,|v_{\boldsymbol iq(i)}|\right)$. Since $\alpha_i\geq0$, this means that either $|v_{\boldsymbol ip(i)}|=1\geq\hat\gamma_0$, or $|v_{\boldsymbol ip(i)}|\geq\hat\gamma_0|v_{\boldsymbol iq(i)}|=\hat\gamma_0$. Thus, $|v_{\boldsymbol i p(i)}|\geq\hat\gamma_0$ in either subcase. Similarly, \eqref{second_thing} leads to the conclusion that $|z_{\boldsymbol i\tilde p(i)}|\geq\hat\gamma_0$.

For the second part, suppose $v_{\boldsymbol ip(i)}<\hat\gamma_0$. Then $z_{\boldsymbol i\tilde p(i)}\geq\hat\gamma_0$ by the first assertion. We now consider two cases. 

Case (a): $\max_{f\in[l]}|v_{\boldsymbol if}|=0$, implying that $\alpha_i=0$. By \eqref{second_thing}, $\zeta_i\neq 0$. Thus, $\zeta_i/\alpha_i=\infty>\hat\gamma_0$.

Case (b): $\max_{f\in[l]}|v_{\boldsymbol if}|>0$. Then $\max_{f\in[l]}|v_{\boldsymbol if}|=1$ due to the Elimination Method. Furthermore, $\hat\gamma_0<1$ by the definitions of $\gamma_0$ and $\hat\gamma_0$. 
Now, the assumption $v_{\boldsymbol ip(i)}<\hat\gamma_0$ and the facts $\max_{f\in[l]}|v_{\boldsymbol if}| =1$, $\hat\gamma_0<1$ and $\max_{f\in[l]}|v_{\boldsymbol if}|=\max(|v_{\boldsymbol ip(i)}|,|v_{\boldsymbol iq(i)}|)$ together imply that $|v_{\boldsymbol iq(i)}|=1$. Consequently, $v_{\boldsymbol ip(i)}<\hat\gamma_0|v_{\boldsymbol iq(i)}|$. Therefore, by Lemma \ref{sign_lemma}, 
\begin{align*}
\zeta_iz_{\boldsymbol i\tilde p(i)}&\geq\frac{\gamma_0|\alpha_iv_{\boldsymbol iq(i)}|-\max\left(\alpha_iv_{\boldsymbol ip(i)},0\right)}{\gamma_0+1}\\
&=\frac{\gamma_0\alpha_i-\max\left(\alpha_iv_{\boldsymbol ip(i)},0\right)}{\gamma_0+1}.
\end{align*}
If $v_{\boldsymbol ip(i)}\leq0$, then the above yields:
\begin{align*}
\frac{\zeta_i}{\alpha_i}&\geq\frac{\gamma_0}{(\gamma_0+1)z_{\boldsymbol i\tilde p(i)}}\geq\frac{\gamma_0}{\gamma_0+1}>\hat\gamma_0
\end{align*}
because $0<z_{\boldsymbol i\tilde p(i)}\leq\max_{f\in[l]}|z_{\boldsymbol if}|=1$. On the other hand, if $v_{\boldsymbol ip(i)}>0$, then
\begin{align*}
\frac{\zeta_i}{\alpha_i}&\geq\frac{\gamma_0-v_{\boldsymbol ip(i)}}{\gamma_0+1}>\frac{\gamma_0-\hat\gamma_0}{\gamma_0+1}=\hat\gamma_0.\end{align*}\end{proof}

The next lemma is the last technical lemma. It forms the crux of the main results of this subsection.

\begin{lemma}\label{lem:lem_9_is_useful}
Let $\{\lambda_i\}_{i=1}^{\tau}$ and $\{U_i\}_{i=1}^\tau$ be such that $|\lambda_i|<1$ and $U_i$ is a linear subspace of $\R^l$ for each $i\in[\tau]$, where $l\in\mathbb N$. Further, suppose $\lambda_i>0$ for some $i\in[\tau]$. Let $\mathcal K$ be the set of all $k_M\in\mathbb N$ such that 
\begin{align} \label{eq:28}
     \sum_{i=1}^{\tau} \lambda_i^k u_{\boldsymbol ij}<\delta\quad\text{for all }1\leq k<k_M\text{ and all }j\in [l],
 \end{align}
 and 
  \begin{align}\label{eq:29}
     \sum_{i=1}^{\tau} \lambda_i^{k_M}u_{\boldsymbol it}\geq \delta\quad\text{ for some }t\in [l],
 \end{align}
 hold for some $\delta\in\R$ and some $(u_{\boldsymbol 1}, \ldots, u_{\boldsymbol {\tau}})\in \prod_{i=1}^{\tau}U_i$.   If $\sup \mathcal K=\infty$, then there exists a $d\in[\tau]$ such that $\lambda_d>0$ and a corresponding non-zero vector $v\in U_d$ such that $v_iv_j\geq 0$ for all $i,j\in [l]$. 
\end{lemma}

\begin{proof} Suppose the lemma is false, i.e., $\sup\mathcal K=\infty$, and $v_iv_j<0$ for some $i,j\in [l]$ whenever there exists a $d\in\tau$ such that $v\in U_d\setminus\{0\}$ and $\lambda_d>0$. The rest of the proof is organized into six steps.

  \textit{Step 1:} By Lemma \ref{lem:gamma_lemma}, there exists a positive constant $\gamma_{d}$ that lower bounds the ratios $|\max_{p\in [l]}v_p|/|\min_{q\in[l]}v_q|$ and $|\min_{p\in [l]}v_p|/|\max_{q\in[l]}v_q|$ for all $v\in U_d$ satisfying $v_{[l]}\neq 0$. Since $\tau<\infty$, the positive constant $\gamma_0:=\min_{d\in [\tau]:\lambda_d>0} \gamma_d$ lower bounds these ratios for every $d\in[\tau]$ for which $\lambda_d>0$. Thus, every non-zero vector lying in $U_d$ has both positive and negative entries that are significant in magnitude.

\textit{Step 2:} It is clear that if $\{\lambda_i\}_{i=1}^\tau$ are distinct, then for $k\gg 1$, $S_j[k]$ will have just one significant term, thereby simplifying our analysis. Since this assumption is invalid, we proceed as follows: pick any $k_M\in\mathcal K\cap [4,\infty)$, let $(u_{\boldsymbol 1},u_{\boldsymbol 2},\ldots,u_{\boldsymbol \tau})$ and $\delta$ be such that \eqref{eq:28} and \eqref{eq:29} hold, and perform the Elimination Method so as to express $S_j[k]$ as:
\begin{align}\label{op_eq}
S_j[k]=\sum_{i=1}^m\alpha_i \mu_i^k v_{\boldsymbol ij}\quad\forall\,\,j\in[l],\,k\leq k_M:k\text{ is even},
\end{align}
and for odd $k\geq 1$, independently as:
\begin{align}
S_j[k]=\sum_{i=1}^m\zeta_i \mu_i^k z_{\boldsymbol ij}\quad\forall\,\,j\in[l],\,k\leq k_M:k\text{ is odd}.
\end{align}
Note that $\mu_i<1$ 
for all $i\in [m]$ because $|\lambda_i|<1$ for $i\in [m]$. Next, note that $\{\alpha_i\}_{i\in [m]}$ and $\{\zeta_i\}_{i\in [m]}$ are determined completely by $(u_{\boldsymbol 1},u_{\boldsymbol 2},\ldots,u_{\boldsymbol \tau})$. Also, w.l.o.g., we assume that $k_M$ is even (otherwise, we can set $u'_{\boldsymbol i}=\lambda_iu_{\boldsymbol i}$ and $k_M'=k_M-1$ so that \eqref{eq:28} and \eqref{eq:29} hold for the primed variables). Since $k_M\geq 4$, \eqref{eq:28} and \eqref{eq:29} imply that $\alpha_i\neq 0$ for some $i\in[m]$ and hence, $\alpha_i>0$ for some $i\in[m]$. All of this implies that for $k\gg1$, if there exist $i\in[\tau]$ and $j\in[l]$ such that $v_{\boldsymbol{i}j}>0$ and $|\alpha_i\mu_i^k v_{\boldsymbol{i}j}|$ is much greater than $\delta$ as well as other terms in $S_j[k]$, then \eqref{eq:28} will be violated..   

\textit{Step 3:} This motivates us to identify an index $s\in [m]$ such that the greatest entry of $v_{\boldsymbol{s}}$, say $v_{\boldsymbol sp}$, is comparable to $1$ and dominates the corresponding sum $S_p[k]$. For this purpose, we wish to ascertain that its weight $\alpha_s$ is comparable to all weights $\alpha_i$ for $i=s+1,\ldots,m$ and  is much greater than the weights $\alpha_i$ for $i=1,\ldots, s-1$ as $\mu_i> \mu_s$ for $i<s$. It is also helpful to compare $v_{s\boldsymbol p}$ with $\sum_{i=s}^m|v_{\boldsymbol ip}|$. With this in mind, we let $\hat\gamma_0:=\gamma_0/(2+\gamma_0)$, and for each $s\in[m]$ that satisfies $\alpha_s> 0$, we define the following quantities: 
\begin{gather*}
\rho_{s1}=\max_{i\leq s-1}\frac{\alpha_i}{\alpha_s},\quad\rho_{s2}=\min_{i\geq s}\frac{\alpha_s}{\alpha_i},\\
p(s)\in \argmax_{j\in [l]}v_{\boldsymbol sj},\quad q(s)\in \argmin_{j\in [l]}v_{\boldsymbol sj},\\
v_{\boldsymbol s0}=
	\begin{cases}
      v_{\boldsymbol sp(s)}, & \text{if}\ v_{\boldsymbol sp(s)}\geq\hat\gamma_0 \\
      1, & \text{otherwise}
	\end{cases},\text{ and}\\
\tau_{s}=\frac{v_{\boldsymbol s0}}{\max_{r=1}^l\left(\sum_{i=s}^m{ |v_{\boldsymbol ir}|}\right)}.
\end{gather*}
Similarly, for each $s\in[m]$ satisfying $\zeta_s\neq 0$, we define:
\begin{gather*}
\tilde\rho_{s1}=\max_{i\leq s-1}\frac{\zeta_i}{\zeta_s},\quad\tilde\rho_{s2}=\min_{i\geq s} \frac{\zeta_s}{\zeta_i},\\
\tilde p(s)\in \argmax_{j\in [l]}z_{\boldsymbol sj},\quad\tilde q(s)\in \argmin_{j\in [l]}z_{\boldsymbol sj},\\
z_{\boldsymbol s0}=
	\begin{cases}
      z_{\boldsymbol s\tilde p(s)}, & \text{if}\ z_{\boldsymbol s\tilde p(s)}\geq\hat\gamma_0 \\
      1, & \text{otherwise}
	\end{cases}, \text{ and}\\
\tilde\tau_{s}=\frac{z_{\boldsymbol s0}}{\max_{r=1}^l\left(\sum_{i=s}^m |z_{\boldsymbol ir}|\right)}.
\end{gather*}
We also let $\rho_{11}=0$ if $\alpha_1>0$ and $\rho_{m2}=1$ if $\alpha_m>0$. Similarly, $\tilde\rho_{11}=0$ if $\zeta_1>0$ and $\tilde\rho_{m2}=1$ if $\zeta_m>0$.


Now, let us bound $\tau_s$:
\begin{align}\label{ineq:tau_bound_1}
\tau_s&\stackrel{(a)}{\geq}\frac{\hat\gamma_0}{\max_{r=1}^l\left(\sum_{i=s}^m |v_{\boldsymbol ir}|\right)}\stackrel{(b)}{\geq}\frac{\hat\gamma_0}{(m-s+1)}\stackrel{}\geq\frac{\hat\gamma_0}{n}.
\end{align}
Here, (a) holds because $\hat\gamma_0\leq v_{\boldsymbol s0}$, and (b) holds because $|v_{j}^{(i)}|\leq1$. On the other hand,
\begin{align}\label{ineq:tau_bound_2}
    \tau_s\leq \frac{1}{\max_{r=1}^l\left(\sum_{i=s}^m |v_{\boldsymbol ir}|\right)}\leq\frac{1}{\max_{j=1}^l|v_{\boldsymbol sj}|}=1.
\end{align}

\textit{Step 4:} We now analyze the evolution of the quantities defined above as $k_M\rightarrow\infty$. Consider any sequence, $\{y^{(h)}\}_{h=1}^\infty=\{(u_{\boldsymbol 1}^{(h)}, \ldots, u_{\boldsymbol {\tau}}^{(h)},\delta^{(h)})\}_{h=1}^\infty$, of variables associated with an increasing and unbounded sequence of solutions $\{k_M^{(h)}\}_{h=1}^\infty\subset \mathcal K$. Since $m<\infty$, there exists an index $M_e\in[m]$ and a subsequence $\{y^{(h_g)}\}_{g=1}^\infty$ of the original sequence $\{y^{(h)}\}_{h=1}^\infty$ such that $M_e\in \argmax_{i\in[m]}\alpha_i^{(h_g)}$ (where $\alpha_i^{(h)}:=\alpha_i(y^{(h)})$), for all $g\in\mathbb N$. Pick such a subsequence and relabel it as $\{y^{(h)}\}_{h=1}^\infty$, so that $0\leq\alpha^{(h)}_i/\alpha^{(h)}_{M_e}\leq1$ for all $i\in[m]$.
 Now that $\{\alpha^{(h)}_i/\alpha^{(h)}_{M_e}\}_{h=1}^\infty$ is bounded for each $i\in[m]$, we may assume (by passing to yet another subsequence if necessary) that $\eta_i:=\lim_{h\rightarrow\infty}\alpha^{(h)}_i/\alpha^{(h)}_{M_e}$ exists for each $i\in [m]$. 

Now, let $r=\min\{i\in[m]:\eta_i>0\}$. Since $\mu_i$ decreases with $i$, we observe that $r$ \textit{indexes the most dominant vector among those that continue to survive even as we increase} $k_M^{(h)}$ (or as we increase $h$). Then $\eta_b=0\text{ for }b\in[i-1]$ and hence, $\lim_{h\rightarrow\infty}\rho^{(h)}_{r1}=0$ and $\lim_{h\rightarrow\infty}\rho^{(h)}_{r2}=\eta_r>0$. Thus, there exists an $h_0\in\mathbb N$ such that for all $h\geq h_0$:
\begin{align}\label{ineq:rho_eta}
    \eta_r/2<\rho_{r2}^{(h)}\leq3\eta_r/2,
\end{align}
 $k^{(h)}_M$ is large enough, and $\rho_{r1}^{(h)}$ is small enough (as will be made precise later). \\

\textit{Step 5:} Our next goal is to show that the greatest positive entry of the dominant vectors is eventually upper bounded by $\hat\gamma_0$. We first restrict ourselves to even-$k$ vectors. We assume $h\geq h_0$, drop the superscript $^{(h)}$ to reduce clutter in notation, and show that:\begin{align}\label{eqn:vpr}
    v_{\boldsymbol rp(r)}<\hat\gamma_0. 
\end{align}



We will assume the contrary and show that $v_{\boldsymbol r}$ dominates other vectors for a range of values of $k$. We will then show that for \eqref{eq:28} to hold, the contribution of $v_{\boldsymbol r}$ should be upper bounded by some function of $\delta$, whereas for \eqref{eq:29} to hold, it should also be lower bounded by a quantity that approaches 0 as $h\rightarrow\infty$. To begin, let $p=p(r)$, suppose $v_{\boldsymbol rp}\geq\hat\gamma_0$ so that $v_{\boldsymbol r0}=v_{\boldsymbol rp}$ and assume that $k$ is even. Then, by \eqref{op_eq}:
\begin{equation}\label{what}
S_j[k]=\sum_{i=1}^{r-1}\alpha_i\mu_i^kv_{\boldsymbol ij}+\alpha_r\mu_r^kv_{\boldsymbol rj}+\sum_{i=r+1}^m\alpha_i\mu_i^kv_{\boldsymbol ij},
\end{equation}
and by \eqref{eq:28}, this implies:
\begin{equation}\label{short0}
\sum_{i=1}^{r-1}\alpha_i\mu_i^kv_{\boldsymbol ij}+\alpha_r\mu_r^kv_{\boldsymbol rj}+\sum_{i=r+1}^m\alpha_i\mu_i^kv_{\boldsymbol ij}<\delta,
\end{equation}
for $k$ in the range $2\leq k< k_M$ and $j\in[l]$. 


Observe that for any $j\in [l]$:
\begin{align}\label{long1}
&\left|\sum_{i=1}^{r-1}\alpha_i\mu_i^kv_{\boldsymbol ij}+\sum_{i=r+1}^m\alpha_i\mu_i^kv_{\boldsymbol ij}\right|\nonumber\\
&\stackrel{(a)}{\leq}\left(\sum_{i=1}^{r-1}\alpha_i|v_{\boldsymbol ij}|\right)\mu_1^k+\left(\sum_{i=r+1}^m\alpha_i|v_{\boldsymbol ij}|\right)\mu_{r+1}^k\nonumber\\
&\stackrel{(b)}{\leq}\left(\sum_{i=1}^{r-1}|v_{\boldsymbol ij}|\right)\rho_{r1}\alpha_r\mu_1^k+\left(\sum_{i=r+1}^m|v_{\boldsymbol ij}|\right)\rho_{r2}^{-1}\alpha_r\mu_{r+1}^k\nonumber\\
&\stackrel{(c)}{\leq}r\rho_{r1}\alpha_r\mu_1^k+\left(\sum_{i=r+1}^m|v_{\boldsymbol ij}|\right)\rho_{r2}^{-1}\alpha_r\mu_{r+1}^k\nonumber\\
&\stackrel{(d)}{\leq}r\rho_{r1}\alpha_r\mu_1^k+(\rho_{r2}\tau_r)^{-1}\alpha_r\mu_{r+1}^kv_{\boldsymbol r0},
\end{align}
where (a) is due to the ordering of the set $\{\mu_i\}_{i=1}^m$, Triangle Inequality and the fact that $\alpha_i\geq0$, (b) follows from the definitions of $\rho_{r1}$ and $\rho_{r2}$, (c) follows from the fact that $\max_{t=1}^l |v_{\boldsymbol it}|=1$, and (d) follows from the definitions of $\tau_r$ and $v_{\boldsymbol r0}$.

Now, we identify a range of $k$ over which the contribution from $\mu_r$ dominates the contributions from both $\mu_{r+1}$ and $\mu_1$. Let $k_{re}:=\max(0,2\lceil0.5\log_{(\mu_r/\mu_{r+1})}(40/\eta_{r}\tau_r)\rceil)$ and $k_r':=2\lfloor0.5\log_{(\mu_1/\mu_{r+1})}(\frac{v_{\boldsymbol r0}}{r\rho_{r1}\rho_{r2}\tau_r})\rfloor$. Then, for $\rho_{r1}$ small enough, \eqref{ineq:tau_bound_1}, \eqref{ineq:tau_bound_2}, and \eqref{ineq:rho_eta} ensure that $k_{re}\leq k_r'< \infty$, and
\begin{equation}\label{short1}
r\rho_{r1}\alpha_r\mu_1^k\leq(\rho_{r2}\tau_r)^{-1}\alpha_r\mu_{r+1}^kv_{\boldsymbol r0}\text{ for }k\leq k_r',
\end{equation}
Furthermore, the definition of $k_{re}$ and \eqref{ineq:rho_eta} imply that
\begin{equation}\label{short2}
(\rho_{r2}\tau_r)^{-1}\alpha_r\mu_{r+1}^k v_{\boldsymbol r0}\leq0.05\alpha_r\mu_r^kv_{\boldsymbol r0}\text{ for }k\geq k_{re}.
\end{equation}
Combining \eqref{long1}, \eqref{short1} and \eqref{short2} yields:
\begin{equation}\label{short3}
\left|\sum_{i=1}^{r-1}\alpha_i\mu_i^kv_{\boldsymbol ij}+\sum_{i=r+1}^m\alpha_i\mu_i^kv_{\boldsymbol ij}\right|\leq0.1\alpha_r\mu_r^kv_{\boldsymbol r0}
\end{equation}
for $k_{re}\leq k\leq k_r'$ and $j\in[l]$. Thus, \textit{if \eqref{eqn:vpr} fails, then the contribution of the dominant vector $v_{\boldsymbol r}$ is much greater than the combined contributions of other even-$k$ vectors when} $k_{re}\leq k\leq k_r'$. Now, \eqref{short3}, the assumption $v_{\boldsymbol rp}\geq\hat\gamma_0$, and  \eqref{short0} at $j=p$ together result in the following:
\begin{align}\label{point_nine}
\delta&>\alpha_r\mu_r^kv_{\boldsymbol rp}-\left |\sum_{i=1}^{r-1}\alpha_i\mu_i^kv_{\boldsymbol ij}+\sum_{i=r+1}^m\alpha_i\mu_i^kv_{\boldsymbol ip} \right |\nonumber\\
&\geq 0.9\alpha_r\mu_r^kv_{
\boldsymbol r0}\text{ for }k_{re}\leq k\leq k_r'. 
\end{align}

By \eqref{what}, \eqref{short3}, and \eqref{point_nine}, we have:
\begin{align}\label{shortening}
S_j[k]&\leq\alpha_r\mu_r^kv_{\boldsymbol rj}+\left |\sum_{i=1}^{r-1}\alpha_i\mu_i^kv_{\boldsymbol ij}+\sum_{i=r+1}^m\alpha_i\mu_i^kv_{\boldsymbol ij}\right|\nonumber\\
&\leq1.1\alpha_r\mu_r^kv_{\boldsymbol r0}\\
&\leq1.1\alpha_r\mu_r^{k_{re}+\log_{1/\mu_r}(11/9)}v_{\boldsymbol r0}<\delta,
\end{align}
for all $j\in[l]$ and $k\in\left[k_{re}+\log_{1/\mu_r}(11/9), \min(k_r',k_M)\right]$ (which is a non-empty interval for a small enough $\rho_{r1}$ and a large enough $k_M$ due to the definition of $k_r'$). 
 
In particular, 
\begin{equation}\label{ineq_limit}
S_j^{(h)}[k_m^{(h)}]\leq 1.1\alpha^{(h)}_r\mu_r^{k_m^{(h)}}v_{\boldsymbol r0}<\delta^{(h)},
\end{equation}
where $k_m^{(h)}:=\min(k_r'^{(h)},k_M^{(h)})$. 
On the other hand, for \eqref{eq:29} to hold for an arbitrarily large $k_M$, we need $v_{\boldsymbol rp}$ to be \textit{much greater} than $\delta$ so as to compensate for the corresponding (small) value of $\mu_r^{k_M}$. This leads to a contradiction. To elaborate, let $t^{(h)}\in [l]$ be the index satisfying \eqref{eq:29}. Then, for every $h\geq h_0$, \eqref{eq:29}, \eqref{ineq_limit}, and \eqref{point_nine} imply:
\begin{align}\label{ineq_clinch}
    \sum_{i=1}^{m}&\alpha_i^{(h)}u_{\boldsymbol it^{(h)}}^{(h)}\left(\mu_i^{k_M^{(h)}}-\mu_i^{k_m^{(h)}}\right)\cr
    &=S_{t^{(h)}}^{(h)}[k_M^{(h)}]-S_{t^{(h)}}^{(h)}[k_m^{(h)}]\cr
    &\geq\delta^{(h)}-1.1\alpha_r^{(h)}\mu_r^{k_m^{(h)}}v_{\boldsymbol r0}\cr
    &> 0.9\alpha_r^{(h)}\mu_r^{k_{re}}v_{\boldsymbol r0}-1.1\alpha_r^{(h)}\mu_r^{k_m^{(h)}}v_{\boldsymbol r0}.
\end{align}
Division by $\alpha_r^{(h)}$ and rearranging the terms yield:
\begin{align}\label{eq:even_clinch}
    1.1\mu_r^{k_m^{(h)}}v_{\boldsymbol r0}+\sum_{i=1}^{m}&\frac{\alpha_i^{(h)}}{\alpha_r^{(h)}}u_{\boldsymbol it^{(h)}}\left(\mu_i^{k_M^{(h)}}-\mu_i^{k_m^{(h)}}\right)>0.9\mu_r^{k_{re}}v_{\boldsymbol r0}.
\end{align}
However, the left-hand side of \eqref{eq:even_clinch} tends to zero as $h\rightarrow\infty$ (since $\eta_r>0$) because $\lim_{h\rightarrow\infty}\rho_{r1}^{(h)}=0$ implies that $\lim_{h\rightarrow\infty}k_r'^{(h)}=\infty$ and in turn that $\lim_{h\rightarrow\infty}k_m^{(h)}=\infty$, whereas the right-hand side remains positive. This contradicts our assumption on $v_{\boldsymbol rp}$, thus proving \eqref{eqn:vpr}.

Now, we establish the odd-$k$ analog of \eqref{eqn:vpr}. Note that the assumption that $k_M$ is even forbids us from repeating our previous arguments. 

By \eqref{ineq:rho_eta}, \eqref{eqn:vpr} and Lemma \ref{lem:sign_extended}, we have
\begin{align} \label{eq_that_clinches_it}
    \zeta_r^{(h)}\geq\hat\gamma_0\alpha_r^{(h)}>0 
\end{align}
for $h\geq h_0$. Therefore, analogous to $M_e$, $\eta_i$ for $i\in[m]$, and $r$, we define $M_o:=\argmax_{i\in[m]}\zeta_{i}^{(h)}$ for $h\geq h_0$, $\tilde\eta_i:=\lim_{h\rightarrow\infty}\zeta^{(h)}_i/\zeta^{(h)}_{M_o}$ for $i\in[m]$, and $\tilde r:=\min\{i\in [m]:\tilde\eta_i>0\}$, respectively (by passing to a subsequence of $\{y^{(h)}[0]\}_{h=1}^\infty$ if necessary). Also, note that we did not use the assumption that $k_M^{(h)}$ is even until \eqref{ineq_limit}. This implies that if  $z_{\tilde{\boldsymbol{r}}\tilde p(\tilde r)}\geq\hat\gamma_0$ holds, then similar to  $\delta^{(h)}> 0.9\alpha_{ r}^{(h)}\mu_{r}^{k_{re}}v_{{\boldsymbol r}0}$, we  have:
\begin{align}\label{eq:analog}
    \delta^{(h)}> 0.9\zeta_{\tilde r}^{(h)}\mu_{\tilde r}^{k_{ro}}z_{\tilde{\boldsymbol r}0},
\end{align}
where $k_{ro}:=\max(0,2\lceil0.5\log_{(\mu_{\tilde r}/\mu_{\tilde r+1})}(40/\tilde\eta_{\tilde r}\tilde\tau_{\tilde r})\rceil)$. On the other hand, $S_{t^{(h)}}[k_M^{(h)}]\geq \delta^{(h)}$  implies:
\begin{align}\label{eq:analog_compare}
    \delta^{(h)}\leq \sum_{i=1}^m \alpha_i^{(h)}\mu_i^{k_M^{(h)}}v_{\boldsymbol ij}\leq \alpha_{M_e}^{(h)}\mu_1^{k_M^{(h)}}m
\end{align}
since $|v_{\boldsymbol ij}|\leq 1$. Then, \eqref{eq:analog} and \eqref{eq:analog_compare} result in:
$$
\frac{\alpha_{M_e}^{(h)}}{\zeta_{\tilde r}^{(h)}}\geq\frac{0.9\mu_r^{k_{ro}}z_{\boldsymbol r0}}{m}\left(\frac{1}{\mu_1}\right)^{k_M^{(h)}},
$$
implying that $\lim_{h\rightarrow\infty} (\alpha_{M_e}^{(h)}/\zeta_{\tilde r}^{(h)})=\infty$. Hence:
\begin{align*}
    \lim_{h\rightarrow\infty}\frac{\zeta_r^{(h)}}{\alpha_r^{(h)}}&=\lim_{h\rightarrow\infty}\frac{\zeta_r^{(h)}}{\zeta_{M_o}^{(h)}}\cdot\frac{\zeta_{M_o}^{(h)}}{\zeta_{\tilde r}^{(h)}}\cdot\frac{\zeta_{\tilde r}^{(h)}}{\alpha_{M_e}^{(h)}}\cdot\frac{\alpha_{M_e}^{(h)}}{\alpha_r^{(h)}}\cr
    &=\tilde\eta_r\cdot\tilde\eta_{\tilde r}^{-1}\cdot0\cdot\eta_r^{-1}=0
\end{align*}
because $\tilde\eta_{\tilde r}$ and $\eta_r$ are positive by the definitions above. However, this would have contradicted \eqref{eq_that_clinches_it}. Therefore:
\begin{align} \label{eqn:zpr}
    z_{\tilde{\boldsymbol{r}}\tilde p(\tilde r)}<\hat\gamma_0.
\end{align}
\textit{Step 6:} Note that \eqref{eqn:vpr}, \eqref{eqn:zpr}, and Lemma \ref{lem:sign_extended} imply that $r\neq \tilde r$. We may assume that $r<\tilde r$ because the case $r>\tilde r$ can be handled similarly. Then by the definition of $\tilde r$, we have $\lim_{h\rightarrow\infty}\zeta_{r}^{(h)}/\zeta_{\tilde r}^{(h)}=0$. Furthermore, by applying Lemma \ref{lem:sign_extended} to both $r$ and $\tilde r$, we obtain $\min(\alpha_{\tilde r}^{(h)}/\zeta_{\tilde r}^{(h)},\zeta_{ r}^{(h)}/\alpha_r^{(h)})\geq\hat\gamma_0$. Therefore, 
\begin{align*}
\lim_{h\rightarrow\infty}\frac{\alpha_{\tilde r}^{(h)}}{\alpha_{M_e}^{(h)}}&=\lim_{h\rightarrow\infty}\frac{\alpha_{\tilde r}^{(h)}}{\zeta_{\tilde r}^{(h)}}\cdot\frac{\zeta_{\tilde r}^{(h)}}{\zeta_{r}^{(h)}}\cdot\frac{\zeta_{r}^{(h)}}{\alpha_{r}^{(h)}}\cdot\frac{\alpha_{r}^{(h)}}{\alpha_{M_e}^{(h)}}\\
&\geq\hat\gamma_0\cdot\infty\cdot\hat\gamma_0\cdot\eta_r=\infty
\end{align*}
because $\eta_r>0$ by our definition of $r$. But this contradicts the definition of $M_e$, thereby proving the lemma. \end{proof}

\begin{figure}[h]
    \centering
        \includegraphics[scale=0.33]{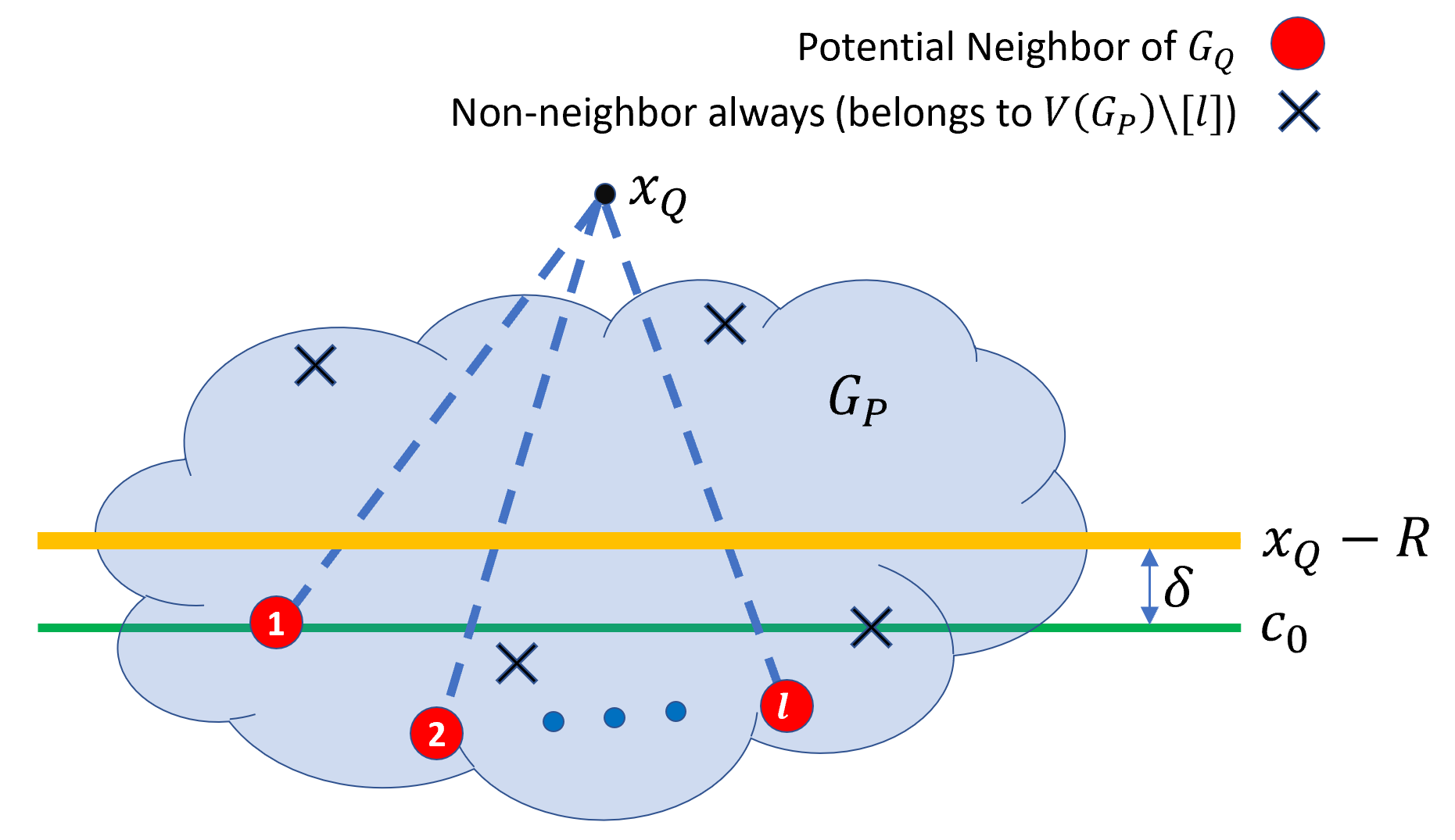}
        \caption{Illustration for the Proof of Lemma \ref{lem:c_below}}
        \label{fig_necessary_conditions}
\end{figure}

We can now state the first main result of this subsection.
\begin{lemma}\label{lem:c_below}
 For every initial state $x[0]\in\R^n$, let $G_P[k]=G_P(x[k])=(V_P, E_P(x[k]))$ and $G_Q[k]=G_Q(x[k])=(V_Q, E_Q(x[k]))$ be two vertex-disjoint induced subgraphs of $\tilde G[k]$ such that $G_{P_0}$, the subgraph of $\Gph$ induced by $V_P$, is connected. Also, let $\mathcal X$ denote the set of all $x[0]\in\mathbb R^{n}$ satisfying assumptions below:
\begin{enumerate}[(a).]
\item \label{item:a}All the agents of $G_Q[0]$ have the same opinion value, i.e., $x_i[0]=x_Q$ for all $i\in V_Q$, where $x_Q\in\mathbb R$ is constant in time but depends on $x[0]$.
\item \label{item:b}$G_P[k]$ and $G_Q[k]$ merge for the first time ever at some time $\kappa(x[0])\in\mathbb N$.
\item \label{connectedness} $G_P[k]$ is a connected graph for $0\leq k< \kappa(x[0])$.
\item \label{item:c}No link break occurs within $G_P[k]$ until time $\kappa(x[0])$.
\end{enumerate}  Furthermore, for some $l\in\mathbb R$, let $[l]$ index the set of nodes of $G_P[0]$ that are adjacent to one or more nodes of $G_Q[0]$ in the graph $\Gph$, as shown in Fig. \ref{fig_necessary_conditions}.

Now, suppose $\sup_{x[0]\in\mathcal X} \kappa(x[0])=\infty$. Then $A_{P_0}$ has an eigenpair $(\lambda, v)$ such that $0<\lambda<1$, $v_i\neq 0$ for some $i\in [l]$, and $v_iv_j\geq 0$ for all $i,j\in[l]$.
\end{lemma}
\begin{proof} 
By Lemma~\ref{lemma:eigenvalue}, $A_{P_0}$ always has an eigenvalue $\lambda\in(0,1)$. So, if the assertion of this lemma is false, then for every eigenpair $(\lambda,v)$ of $A_{P_0}$ with $0<\lambda<1$, we have $v_pv_q<0$ for some $p,q\in[l]$, while $\sup_{x[0]\in\mathcal X} \kappa (x[0])=\infty$. 


Now, for $G_P$ and $G_Q$ to merge for the first time at $\kappa=\kappa(x[0])$, we require $G_P[k]$ and $G_Q[k]$ to be (I): disconnected from each other in $\tilde G[k]$ until time $\kappa$, and (II): connected with each other in $\tilde G[\kappa]$. 

Given that $G_Q[0]$ is at the consensus state $x_Q$ and the set of potential neighbors of $G_Q$ in $G_P$ is $[l]$, condition (I) is equivalent to:  
\begin{equation}\label{sixteen}
|x_i[k]-x_Q|>R\quad\forall\,\,i\in[l]\text,\,\,k\leq \kappa-1.
\end{equation}

 Since no link break occurs within $G_P[k]$ until the merging event of interest takes place, $G_P[k]$ remains connected in $\tilde G[k]$ until $\kappa(x[0])$. Moreover, since $\sup_{x[0]\in\mathcal X}\kappa(x[0])=\infty$, and because all the intra-component link formations taking place in $\tilde G[k]$ occur in $O(n^8\log n)$ steps as per Proposition \ref{merging_must}, we may choose an $x[0]\in\R^n$ such that $\kappa(x[0])$ is large enough and $G_P[k]$ remains constant and connected during a time interval $[k_c,\kappa(x[0]))$ for some $k_c<\kappa(x[0])$. As a result, we may further assume that for a sufficiently large $\kappa(x[0])$, the sub-network of $\Gph$ corresponding to $G_{P_0}$ achieves $R/4$-convergence to a consensus state $c_0\allone_P$ at some time $k_R\in [k_c,\kappa(x[0]))$. We now shift the origin of our time axis to $k_R$, thus obtaining $G_P[0]=G_{P_0}$. Also, w.l.o.g., we assume $c_0\leq x_Q$. Then \eqref{sixteen} along with $R/4$-convergence together  yield the following necessary condition for (a): 
\begin{equation}\label{late_merge}
x_j[k]<x_Q-R\quad \forall\,\,j\in[l],\,\,k\in\{0, 1, \ldots, k_M-1\},
\end{equation}
where $k_M:=\kappa-k_R$. In this setting, condition (II) is equivalent to $|x_i[k_M]-x_Q|\leq R$ for some $i\in[l]$. An implication is: 
\begin{equation}\label{merge_now}
x_t[k_M]\geq x_Q-R\text{ for some }t\in[l].
\end{equation}
By Lemma \ref{lemma:diagonalizability}, we can express $x[k]$ in terms of the eigenpairs $\{(\lambda_i,u_{\boldsymbol i})\}_{i=1}^{n_P}$ of $A_{P_0}$ in order to rewrite \eqref{late_merge} and \eqref{merge_now} as:
\begin{align} \label{eq:essence_1}
     \sum_{i=2}^{n_P} \lambda_i^k u_{\boldsymbol ij}<\delta\quad\text{for all }0\leq k<k_M\text{ and all }j\in [l],
 \end{align}
 and 
  \begin{align}\label{eq:essence_2}
     \sum_{i=2}^{n_P} \lambda_i^{k_M}u_{\boldsymbol it}\geq \delta\quad\text{ for some }t\in [l],
 \end{align}
where $\delta:=x_Q-R-c_0$ and the sum index $i\geq 2$ because $\lambda_1=1$ and $u_{\boldsymbol 1}=c_0\allone_P$. Since $|\lambda_i|<1$ for all $2\leq i\leq n_P$ and $\max_{i=2}^{n_P}\lambda_i>0$ by Lemma \ref{lemma:eigenvalue}, an application of Lemma \ref{lem:lem_9_is_useful} immediately yields the required condition on $A_{P_0}$.

\end{proof}

We now generalize Lemma \ref{lem:c_below} by allowing both the subgraphs $G_P$ and $G_Q$ to have any of the initial states that force them to remain connected components of the influence graph until they merge (or forever if they do not merge). However, we will need a definition and some notation.

Suppose $G_{P_0}=(V_P,E_{P_0})$ and $G_{Q_0}=(V_Q,E_{Q_0})$ are two induced subgraphs of $\Gph$ such that $V_P\cap V_Q=\emptyset$. Let $\{(i_e,j_e)\}_{e=1}^b\subset V_P\times V_Q$ be the set of boundary edges of $\{G_{P_0},G_{Q_0}\}$ in $\Gph$ (i.e., the set of edges connecting $G_{P_0}$ with $G_{Q_0}$ in $\Gph$), and let $\{1\}\cup\{\lambda_d\}_{d=1}^m$ be the union of the sets of eigenvalues of $A_{P_0}$ and $A_{Q_0}$ (such that $\lambda_d\neq 1$ for all $d\in [m]$). Further, for each $d\in [m]$, let $U_d(P)$ (respectively, $U_d(Q)$) be the eigenspace of $\lambda_d$ with respect to $A_{P_0}$ (respectively, $A_{Q_0}$) if $\lambda_d$ is an eigenvalue of $A_{P_0}$ (respectively, $A_{Q_0}$), and let $U_d(P)=\{0\}$ (respectively, $U_d(Q)=\{0\}$), otherwise. Finally, for each $d\in [m]$, let $f_{br}^P(u):=[u_{i_1}\,\,\ldots\,\,u_{i_b}]^T$ for all $u\in U_d(P)$, and let $f_{br}^Q(w):=[w_{j_1}\,\,\ldots\,\,w_{j_b}]^T$ for all $w\in U_d(Q)$. Note that the dimensions of $f_{br}^P(u)$ and $f_{br}^Q(w)$ equal $b$ for all $u\in U_d(P)$ and $w\in U_d(Q)$.

\begin{definition} For each $d\in [m]$, the \textit{boundary-restricted eigenspace} of $\lambda_d$ associated with $\{G_{P_0},G_{Q_0}\}$ is the set $\hat U_d^{PQ}:= \hat U_d^P+ \hat U_d^Q$, where $\hat U_d^P := \text{span}(\{f_{br}^P(v):v\in U_d(P)\})$ and $\hat U_d^Q := \text{span}(\{f_{br}^Q(v):v\in U_d(Q)\})$. We refer to any vector $v\in \hat U_d^{PQ}$ as a \textit{boundary-restricted eigenvector} of $\{G_{P_0},G_{Q_0}\}$ corresponding to the eigenvalue $\lambda_d$.
\end{definition} 

\begin{proposition}
\label{prop:final}
For every initial state $x[0]\in\R^n$, let $G_P[k]=G_P(x[k])=(V_P, E_P(x[k]))$ and $G_Q[k]=G_Q(x[k])=(V_Q, E_Q(x[k]))$ be two vertex-disjoint induced subgraphs of $\tilde G[k]$, and let $\mathcal X$ denote the set of all $x[0]\in\mathbb R^n$ satisfying the assumptions below:
\begin{enumerate}[(i).]

\item \label{item:i}$G_P[k]$ and $G_Q[k]$ merge at time $k_M(x[0])\in\mathbb N$ for the first time.
\item \label{connectedness_2} $G_P[k]$ and $G_Q[k]$ are connected graphs for $0\leq k< k_M(x[0])$.
\item \label{item:ii}No link breaks within $G_P[k]$ or $G_Q[k]$ until time $k_M(x[0])$.
\end{enumerate} Next, let $G_{P_0}$ and $G_{Q_0}$ be the subgraphs of $\Gph$ induced by $V_P$ and $V_Q$, respectively, and let $b$ be the number of boundary edges of $\{G_{P_0},G_{Q_0}\}$ in $\Gph$.
Furthermore, let $\{\lambda_d\}_{d=1}^{m}\cup\{1\}$ be the union of the sets of eigenvalues of $A_{P_0}$ and $A_{Q_0}$ such that $\lambda_d\neq 1$ for $d\in[m]$.

Now, suppose $\sup_{x[0]\in\mathcal X} k_M(x[0])=\infty$. Then there exists an index $d\in[m]$ such that $0<\lambda_d<1$ and a corresponding vector $\hat v\in\hat U_d^{PQ}$ satisfying $\hat v_e\neq 0$ for some $e\in[b]$ and $\hat v_{e}\hat v_{f}\geq0$ for all $e,f\in[b]$.
\end{proposition}

\begin{proof}
Since no link break occurs within $G_P[k]$ or $G_Q[k]$ until they merge, both of them remain connected in $\tilde G[k]$ until $k_M(x[0])$. Moreover, since $\sup_{x[0]\in\mathcal X}k_M(x[0])=\infty$, and because all the intra-component link formations taking place in $\tilde G[k]$ occur in $O(n^8\log n)$ steps as per Proposition \ref{merging_must}, we may choose an $x[0]\in\R^n$ such that $k_M(x[0])$ is large enough, and $G_P[k]$ and $G_Q[k]$ both remain constant and connected during a time interval $[k_c,k_M(x[0]))$ for some $k_c<k_M(x[0])$. As a result, we may further assume that the sub-networks of $\Gph$ corresponding to $G_{P_0}$ and $G_{Q_0}$ achieve $R/4$-convergence to their respective consensus states at some time $k_R\in [k_c,k_M(x[0]))$. We now shift the origin of our time axis to $k_R$, thus obtaining $G_P[0]=G_{P_0}$ and $G_Q[0]=G_{Q_0}$.

Now, we express the initial states of $G_P$ and $G_Q$ as
\begin{align*}
    x_P[0] &= c_P\allone_P + \sum_{d=1}^m\beta_d^Pv_{\boldsymbol d}^P\cr 
    x_Q[0] &= c_Q\allone_Q + \sum_{d=1}^m\beta_d^Qv_{\boldsymbol d}^Q,
\end{align*}
where $c_P,c_Q\in\R$ depend on our choice of $x_P[0]$ and $x_Q[0]$, and the vectors are chosen such that for each $d\in[m]$, $v_{\boldsymbol d}^P$ (respectively, $v_{\boldsymbol d}^Q$) is an eigenvector of $A_{P_0}$ (respectively, $A_{Q_0}$) corresponding to $\lambda_d$ iff $\lambda_d$ is an eigenvalue of $A_{P_0}$ (respectively, $A_{Q_0}$) and $v_{\boldsymbol d}^P=0$ (respectively, $v_{\boldsymbol d}^Q=0$) otherwise. This is possible because $A_{P_0}$ and $A_{Q_0}$ are diagonalizable by Lemma \ref{lemma:diagonalizability}. In addition, we assume that $\{v_{\boldsymbol d}^P\}_{d=1}^{m}\setminus\{0\}$ and $\{v_{\boldsymbol d}^Q\}_{d=1}^m\setminus\{0\}$ are bases of eigenvectors for $A_{P_0}$ and $A_{Q_0}$, respectively. 

Next, let $\{(i_e,j_e)\}_{e=1}^b\subset V_P\times V_Q$ enumerate the set of boundary edges of $\{G_{P_0},G_{Q_0}\}$ in $\Gph$. Note that assumption \eqref{item:i} requires $|x_{i_e}[k]-x_{j_e}[k]|>R$ for all $e\in [b]$ and $0\leq k< k_M:=k_M(x[0])$, and $|x_{i_t}[k_M]-x_{j_t}[k_M]|\leq R$ for some $t\in [b]$. Now, for a given $e\in [b]$, we could either have 
\begin{align}\label{eq:former}
    x_{i_e}[k] - x_{j_e}[k]&>R, \mbox{ or}\\ \label{eq:latter}
    x_{j_e}[k] - x_{i_e}[k]&>R
\end{align}
 for a particular $k\in[0,k_M)$. Suppose \eqref{eq:former} holds at some $k_1\in [0,k_M)$ and \eqref{eq:latter} at some $k_2\in [0, k_M)$. Then $\max(x_{i_e}[k_1] - x_{i_e}[k_2],x_{j_e}[k_2] - x_{j_e}[k_1])>R$. But this contradicts the assumption that both $G_P$ and $G_Q$ have achieved $R/4$-convergence to their respective consensus states at time 0. Therefore, for a given $e\in [b]$, if \eqref{eq:former} holds for some $k\in [0,k_M)$, then it must hold for all $k\in [0,k_M)$. Similarly, we can show that for a given $k\in [0,k_M)$, if \eqref{eq:former} holds for some $e\in [b]$, then it must hold for all $e \in [b]$. The same applies to \eqref{eq:latter}. Hence, w.l.o.g., we assume \eqref{eq:latter} for all $e\in [b]$ and all $0\leq k<k_M$.
 
 Now, for each $d\in [m]$, let $\hat v_{\boldsymbol d}:=[\hat v_{\boldsymbol d1}\,\,\ldots\,\,\hat v_{\boldsymbol db}]^T$, where $\hat v_{\boldsymbol de} = \beta_d^P v_{\boldsymbol di_e} - \beta_d^Q v_{\boldsymbol dj_e}$ for $e\in [b]$. Further, let $\delta:=c_Q-R-c_P$. With these definitions and the assumption given by \eqref{eq:latter}, we can express assumption \eqref{item:i} of the proposition as:
 \begin{align*}
     \sum_{d=1}^m \lambda_d^k\hat v_{\boldsymbol de}<\delta\quad\text{for all }0\leq k<k_M\text{ and all }e\in [b],
 \end{align*}
 and $\sum_{d=1}^m \lambda_d^{k_M}\hat v_{\boldsymbol dt}\geq \delta\quad\text{ for some }t\in [b]$. 
Since $|\lambda_d|<1$ for all $d\in [m]$ and $\max_d\lambda_d>0$ by Lemma \ref{lemma:eigenvalue}, and since $(\hat v_{\boldsymbol 1},\ldots, \hat v_{\boldsymbol{m}})\in\prod_{d=1}^m \hat U_d^{PQ}$, the assertion of Proposition \ref{prop:final} now follows immediately from Lemma \ref{lem:lem_9_is_useful}.
\end{proof}

%% file: cdc_part.tex
\subsection{Graphs with Finite Maximum $\epsilon$-Convergence Time}

We now show that the $\epsilon$-convergence time of a complete $r$-partite graph is bounded. For this purpose, we characterize the eigenvectors of the normalized adjacency matrix of a complete $r$-partite graph that has all the self-loops. 

 For $n\in\mathbb N$, we define a complete $r$-partite graph $G=([n],E)$ to be a graph with partitioning of its vertices into $V_1,\ldots, V_r\subset[n]$, and $(i,j)\in E$ iff $(i,j)\notin\cup_{l=1}^r V_l^2$. Let $G$ have all the $n$ self-loops, let $A\in\R^{n\times n}$ be the normalized adjacency matrix of $G$, and let $n_i:=|V_i|\geq 1$ for $i\in [r]$. For each $i\in[r]$, let $V_i=\{N_{i-1}+1,\ldots,N_{i}\}$, where $N_j:=\sum_{i=1}^j n_i$ for $j\in [r]$ and $N_0:=0$. Finally, we define the matrix $B\in\R^{r\times r}$ by: $$B_{ij}:=
    \begin{cases}
        \frac{1}{n-n_i+1} & \text{ if }j=i\\
        \frac{n_j}{n-n_i+1} & \text{ if }j\neq i\\
    \end{cases},
 $$
and let $\{w^{(i)}\}_{i=1}^q$ be an eigenvector basis for $B$ with $\{\lambda^{(i)}\}_{i=1}^q$being the corresponding eigenvalues.
\begin{lemma}\label{lem:evals} The matrices $A$ and $B$ (as described above) have the following properties:
\begin{enumerate}[(i)]
    \item For each $i\in [r]$ such that $n_i\geq 2$ and each $ t\in\{2,\ldots,n_i\}$, the vector $v^{(i,t)}\in\R^n$, defined as:
$$
    v_j^{(i,t)} := 
    \begin{cases}
        +1, & \text{if }j=N_{i-1}+1\\
        -1, & \text{if }j=N_{i-1}+t\\
        0, & \text{otherwise}
    \end{cases},
$$
   is an eigenvector of $A$ corresponding to $1/(n-n_i+1)$. Moreover, the set $U_1:=\{v^{(i,t)}:2\leq t\leq n_i, i\in [r]\}$ is a set of linearly independent vectors.
    \item  For each $i\in [q]$, the vector $\tilde v^{(i)}\in\R^n$, defined as $\Tilde{v}_p^{(i)}=w_j^{(i)}$ for all $p\in V_j$ and $j\in [r]$, is an eigenvector of $A$ corresponding to $\lambda^{(i)}$.
    \item The eigenvectors of $B$ span $\R^r$, i.e., $q=r$.
    \item If $\lambda^{(i)}\neq 1$, then $\lambda^{(i)}\leq0$ for all $i\in[r]$. 
    \item $U:=\cup_{j=1}^r\{v^{(j,t)}:2\leq t\leq n_j\}\cup\{\tilde v^{(i)}\}_{i=1}^r$ is an eigenvector basis for $A$.
\end{enumerate}
\end{lemma}

\begin{proof}
 Observe that for all $p\in[r]$, the degree of each vertex in $V_p$, with its self-loop counted, is $n-n_p+1$. Hence, given $i\in[r]$, for all $p\in [r]\setminus\{i\}$ and $j\in V_p$, we have:
\begin{align*}
    (Av^{(i,t)})_j
    &=\frac{1}{n-n_p+1}\left(v^{(i,t)}_{N_{i-1}+1}+v^{(i,t)}_{N_{i-1}+t}\right)\\
    &=0=\frac{1}{n-n_p+1}v_j^{(i,t)}.
\end{align*}
Next, if $j=N_{i-1}+1$, then
\begin{align*}
    (Av^{(i,t)})_j&=\frac{1}{n-n_i+1}v^{(i,t)}_{N_{i-1}+1}+0\cdot v^{(i,t)}_{N_{i-1}+t} \\
    &=\frac{1}{n-n_i+1}v_j^{(i,t)}.
\end{align*}
Similarly, $(Av^{(i,t)})_j
    =v_j^{(i,t)}/(n-n_i+1)$ also holds for $j=N_{i-1}+t$.
Finally, for $j\in V_i\setminus\{N_{i-1}+1,N_{i-1}+t\}$, we have $(Av^{(i,t)})_j=A_{jj}\cdot 0+\sum_{s\in V\setminus V_i} A_{js}\cdot0=v_j^{(i,t)}/(n-n_i+1)$. So, for each $i\in [r]$ and each $t\in\{2,\ldots,n_i\}$, $v^{(i,t)}$ is an eigenvector of $A$ corresponding to $\frac{1}{n-n_i+1}$. By taking linear combinations, we can easily see that $\{v^{(i,t)}:2\leq t\leq n_i,i\in [r]\}$ are linearly independent vectors. This proves (i). 

 As for (ii), for any $j\in [r]$ and $p\in V_j$, we have:
\begin{align*}
    (&A\tilde v^{(i)})_p=\sum_{s=1}^n
    A_{ps}\tilde{v}^{(i)}_s\cr
    &=\frac{1}{n-n_j+1}\cdot \tilde v^{(i)}_p + \sum_{l\in [r]\setminus\{j\}}\left(\sum_{m\in V_l}\frac{1}{n-n_j+1}\cdot \tilde v^{(i)}_m\right)\cr
    &=\frac{1}{n-n_j+1}\cdot w^{(i)}_j + \sum_{l\in [r]\setminus\{j\}}\frac{n_l}{n-n_j+1}\cdot w^{(i)}_l\cr
    &=\sum_{l=1}^r B_{jl}w^{(i)}_l=(Bw^{(i)})_j=(\lambda^{(i)} w^{(i)})_j = \lambda^{(i)}\tilde v_p^{(i)}.
\end{align*}

 In order to prove (iii), note that $B=D_1SD_2$, where $D_1:=\text{diag}(\frac{1}{n-n_1+1},\ldots,\frac{1}{n-n_r+1})$, $D_2:=\text{diag}(n_1,\ldots,n_r)$, and $S$ is the symmetric $r\times r$ matrix given by:
\begin{align*}
    S_{ij}=
    \begin{cases}
        \frac{1}{n_i} & \text{ if }j=i\\
        1 & \text{ if }j\neq i\\
    \end{cases}.
\end{align*}
Now, observe that the commutativity of diagonal matrices allows us to express $D_1SD_2$ as $D_A(D_BSD_B)D_A^{-1}$, where $D_A:=(D_1D_2^{-1})^{\frac{1}{2}}$ and $D_B:=(D_1D_2)^{\frac{1}{2}}$. Thus, $B=D_A(D_BSD_B)D_A^{-1}$ is similar to the symmetric matrix $D_BSD_B$ and hence, its eigenvectors span $\R^r$, i.e., $q=r$.

As for (iv), for any $\lambda^{(i)}\neq 1$, we know that $D^{\frac{1}{2}}\tilde v^{(i)}$ is an eigenvector of $D^{\frac{1}{2}}AD^{-\frac{1}{2}}$ which is a symmetric matrix as per Lemma 1 of~\cite{parasnis2018hegselmann}. Hence, $\{D^{\frac{1}{2}}\tilde v^{(i)}\}_{i=1}^r$ is an orthogonal set. Since $\allone\in\{\tilde v^{(i)}\}_{i=1}^r$, this implies that \begin{align}\label{eq:pseudo_orthogonality}
    \allone^T D\tilde v^{(i)}=0\text{ if }i\in[r]\text{ and }\lambda^{(i)}\neq 1,
\end{align}thereby forcing each $\tilde v^{(i)}$ to have both positive and negative entries. Now, pick any $i\in[r]$ for which $\lambda^{(i)}\neq 1$, and let $s\in[r]$ be the index such that $w_j^{(i)}\geq0$ for $j\in[s]$ and $w_j^{(i)}<0$ otherwise (we can always label the vertices suitably so that such an $s$ exists). Then \eqref{eq:pseudo_orthogonality} implies that $1\leq s\leq r-1$. Consequently, we have the following relations:
\begin{align*}
    \lambda^{(i)}|w_1^{(i)}|&=\frac{|w_1^{(i)}|+\sum_{j=2}^sn_j|w_j^{(i)}|-\sum_{j=s+1}^rn_j|w_j^{(i)}|}{n-n_1+1}\cr
    -\lambda^{(i)}|w_{s+1}^{(i)}|&=\frac{\sum_{j=1}^sn_j|w_j^{(i)}|-\sum_{j=s+2}^rn_j|w_j^{(i)}|-|w_{s+1}^{(i)}|}{n-n_{s+1}+1}.
\end{align*}
On the basis of this, we have the following for $\lambda^{(i)}\notin \{0,1\}$:
\begin{align*}
    0&< (n-n_1+1)|w_1^{(i)}| + (n-n_{s+1}+1)|w_{s+1}^{(i)}|\cr
    &=-\frac{(n_1-1)|w_1^{(i)}|+(n_{s+1}-1)|w_{s+1}^{(i)}|}{\lambda^{(i)}},
\end{align*}
implying $\lambda^{(i)}<0$ because $n_1,n_{s+1}\geq 1$ by assumption.

For part (v), note that $U_1$ and $\{\tilde v^{(i)}\}_{i=1}^r$ are linearly independent sets by assertions (i) and (ii). Also, observe that
$\text{span}\{\tilde v^{(i)}\mid i\in[r]\}=\text{span}^{\perp}\{v^{(j,t)}\mid j\in[r],t\in \{2,\ldots,n_j\} \}$
because $
    {\tilde {v}^{(i)T}} v^{(j,t)} =w_j^{(i)}\times 1+w_j^{(i)}\times-1=0.$
Finally, noting that $|U|=\sum_{j=1}^r(n_j-1)+r=\sum_{j=1}^rn_j=n$, we conclude that $U$ is an eigenvector basis for $A$. 
\end{proof}

\begin{remark}\label{rem:lemma_1_use}
Points (1), (4) and (5) of Lemma \ref{lem:evals}, along with the fact that eigenspaces are linear, imply that every eigenpair ($\lambda,v)$ of $A$ that satisfies $\lambda\in(0,1)$, corresponds to some $i\in [r]$ such that $|V_i|\geq 2$ and $v_s=0$ for all $s\notin V_i$. Furthermore, $\sum_{s\in V_i}v_s=0$ for such an $i$.
\end{remark}

We are now well equipped to establish our main result.

\begin{proposition} \label{prop:complete_r_partite}
Let $n\in\mathbb N$ and $\epsilon>0$ be given, and let $\Gph=([n],\Eph)$ be a complete $r$-partite graph for some $r\in [n]$. Then $k^*_\epsilon(\Gph)<\infty$. 
\end{proposition}

\begin{proof}
If $\Gph$ is a complete $1$-partite graph, then $\Eph=\emptyset$ and hence $k_\epsilon^*(\Gph)=0$. On the other hand, if $r=n$, then $\Gph=K_n$. In this case, $k_\epsilon^*(\Gph)=O(n^3)<\infty$ by~\cite{mohajer2013convergence} and~\cite{parasnis2018hegselmann}. Therefore, we assume $1<r<n$ hereafter.

Suppose $k_\epsilon(\Gph)=\infty$. From Proposition \ref{merging_must}, we know that  arbitrarily slow convergence happens only in the presence of arbitrarily slow merging and that all the other structural changes in $\tilde G[k]$ occur in $O(n^8\log n)$ steps. Hence, it suffices to show that no two connected components of the influence graph can take an arbitrarily long period of time to merge, under the assumption that no link breaks occur.

For this purpose, let $V_1,\ldots,V_r$ be the $r$ parts of $\Gph$, and let $V_P, V_Q\subset [n]$ be any two disjoint sets. Further, let $\{(i_e,j_e)\}_{e=1}^b\subset V_P\times V_Q$ be the set of boundary edges connecting $G_{P_0}$ and $G_{Q_0}$ in $\Gph$, and let $\{1\}\cup\{\lambda_d\}_{d=1}^m$ be the union of the sets of eigenvalues of $A_{P_0}$ and $A_{Q_0}$ (such that $\lambda_d\neq 1$ for all $d$). Now, since $\Gph$ is a complete $r$-partite graph, it follows that $G_{P_0}$ and $G_{Q_0}$ are also complete $p$-partite and $q$-partite graphs for some $p,q\in [r]$, and their parts are given by the partitions $\{V_P\cap V_i\}_{i=1}^r\setminus\{\emptyset\}$ and $\{V_Q\cap V_i\}_{i=1}^r\setminus\{\emptyset\}$, respectively.

Next, for each initial state $x[0]\in\R^n$, let $G_P[k]=G_P(x[k])=(V_P,E_P(x[k]))$ and $G_Q[k]=G_Q(x[k])=(V_Q,E_Q(x[k]))$ be disconnected from each other in $\tilde G[k]=\tilde G(x[k])$ until they merge at time $k_M(x[0],V_P,V_Q)$. As per our earlier reasoning, we may restrict our attention to the subset $\mathcal X(V_P,V_Q)\subset\R^n$ of initial states for which (i) $k_M(x[0],V_P,V_Q)<\infty$, i.e., merging occurs, (ii) no link breaks occur within $G_P[k]$ or $G_Q[k]$ until they merge, i.e., for $k\leq k_M(x[0], V_P, V_Q)$, and (iii) both $G_P[k]$ and $G_Q[k]$ are connected graphs for $k\leq k_M(x[0],V_P,V_Q)$.

Now, suppose $\sup_{x[0]\in\mathcal{X}{(V_P,V_Q)}}k_M(x[0],V_P,V_Q)=\infty$. Then Proposition \ref{prop:final} implies that there exists a $d\in [m]$ with $\lambda_d\in(0,1)$ and a corresponding vector $v\in \hat U_d^{PQ}$ satisfying $v_e\neq 0$ for some $e\in [b]$ and $v_fv_g\geq 0$ for all $f,g\in [b]$. Since $v_e = u_{i_e}+ w_{j_e}$ for some $u\in U_d(P)$ and $w\in U_d(Q)$, we have either $u_{i_e}\neq 0$ or $w_{j_e}\neq 0$. W.l.o.g., we assume $u_{i_e}>0$ (and hence that $(\lambda_d,u)$ is an eigenpair of $A_{P_0}$). Now, let $\rho\in [r]$ and $\sigma\in [r]$ denote the indices for which $i_e\in V_{P\rho}:=V_P\cap V_\rho$ and $j_e\in V_{Q\sigma}:=V_Q\cap V_{\sigma}$. Then observe that $\rho\neq\sigma$ because $(i_e,j_e)\in \Eph$. Also, by Remark \ref{rem:lemma_1_use}, $\lambda_d\in(0,1)$ implies that $\sum_{s\in V_{P\rho}}{ u_s}=0$. Hence, there exists another node $z\in V_{P\rho}$ such that $u_{z}<0$. Now, two cases arise: either $|V_{Q\sigma}|=1$ or $|V_{Q\sigma}|\geq 2$. 

Consider Case 1: $|V_{Q\sigma}|=1$, i.e., $V_{Q\sigma}=\{j_e\}$. Now, if $\lambda_d$ is not an eigenvalue of $A_{Q_0}$, then $U_d=\{0\}$, which means $w=0$. Hence, $w_{j_e}=0$. Otherwise, by Remark \ref{rem:lemma_1_use}, Lemma \ref{lem:evals} requires $w_{j_e}=0$ because $\lambda_d>0$ and $|V_{Q\sigma}|<2$. Thus, $w_{j_e}=0$ is true whenever $|V_{Q\sigma}|=1$. Moreover, $\rho\neq\sigma$ implies that $(z,j_e)\in \Eph$. Since $z\in V_P$ and $j_e\in V_Q$, we may denote $z$ by $i_f$ and $j_e$ by $j_f$ so that $(z,j_e)$ is the $f$-th boundary edge, $(i_f, j_f)$, for some $f\in [b]$. But now, $v_f=u_{i_f}+w_{j_f}=u_z+w_{j_e}=u_z<0$, whereas $v_e=u_{i_e}>0$. As a result, $v_e v_f<0$, thus contradicting the requirement $v_fv_g\geq 0$ for all $f,g\in [b]$.

On the other hand, in Case 2: $|V_{Q\sigma}|\geq 2$, both $w_{j_e}=0$ and $w_{j_e}\neq 0$ are possible subcases. If $w_{j_e}=0$, then we simply repeat the arguments of the previous paragraph to show that $v_ev_f<0$ for some $f\in [b]$. So, assume $w_{j_e}\neq 0$. Then $(\lambda_d,w)$ is necessarily an eigenpair of $A_{Q_0}$. Therefore, the requirement $\sum_{s\in V_{Q\sigma}}w_s=0$ of Lemma \ref{lem:evals} implies $w_{y}w_{j_e}<0$ for some $y\in V_{Q\sigma}$. First, suppose $w_{j_e}>0$ and $w_y<0$. Then, $\rho\neq\sigma$ implies that $(z,y)\in \Eph$ and hence that $(z,y)$ is a boundary edge. By denoting $(z,y)$ as the $f$-th boundary edge $(i_f,j_f)$ for some $f\in [b]$, we have $v_{f}=u_{i_f}+w_{j_f}=u_z+w_y<0$. However, we still have $v_{e} = u_{i_e}+w_{j_e}>0$, implying that $v_ev_f<0$. Now, assume $w_{j_e}<0$ and $w_y>0$. Then, by denoting the boundary edges $(z,j_e)$ and $(i_e,y)$ as $(i_\alpha, j_\alpha)$ and $(i_\beta,j_\beta)$, respectively for some $\alpha,\beta\in [b]$, we have $v_{\alpha} = u_{z} + w_{j_e}<0$ and $v_\beta = u_{i_e} + w_y>0$. This implies that $v_\alpha v_\beta<0$. Thus, the requirement $v_fv_g\geq 0$ for all $f,g\in [b]$ is violated in Case 2 as well.

Hence, $\sup_{x[0]\in \mathcal X(V_P,V_Q)}k_M(x[0],V_P,V_Q)<\infty$. Note that this applies to every selection of $V_P\subset V$ and $V_Q\subset V$ such that $V_P\cap V_Q=\emptyset$. Moreover, since the number of such choices of $V_P$ and $V_Q$ is finite, we conclude that no merging event can be delayed indefinitely in the social HK dynamics on the given $\Gph$. This completes the proof. 
\end{proof}

%% file: conclusion.tex
\section{CONCLUSION AND FUTURE DIRECTIONS}\label{sec:conclusion}
In this paper, we have investigated the convergence properties of the social HK model of opinion dynamics. We have shown that for certain physical connectivity graphs, we cannot even guarantee $\epsilon$-convergence to the steady state within a bounded time-frame, much less termination in finite time. In addition, we have shown that complete $r$-partite graphs have bounded $\epsilon$-convergence times. Moreover, we can observe that the necessary and sufficient conditions provided by Proposition \ref{prop:positive_evector} and Lemma \ref{lem:c_below} are nearly tight (i.e., tight under the assumption $v_iv_j\neq0$, in addition to the other assumptions made by these two results). However, finding a set of necessary and sufficient conditions for arbitrarily slow merging (and thereby for arbitrarily slow $\epsilon$-convergence) that are tight in the most general case, remains an interesting open problem. Also open is the problem of finding other classes of graphs that have bounded $\epsilon$-convergence times.

